\documentclass[onefignum,onetabnum]{siamart190516}

\usepackage{lmodern}
\usepackage[T1]{fontenc}
\DeclareSymbolFont{myletters}{OML}{ztmcm}{m}{it}
\DeclareMathSymbol{\uplambda}{\mathord}{myletters}{"15}
\def\A{\mathcal{A}}
\def\B{\mathcal{B}}
\def\f{{\scriptstyle{\mathcal{C}}}}
\def\calo{\mathcal{O}}
\def\T{\mathcal{T}}
\def\Dh{D_{h_\T}}
 \usepackage{multirow}
\usepackage[percent]{overpic}
\usepackage{graphicx,graphics,epsfig}
\usepackage{amsmath,amssymb,dsfont}
\usepackage{mathrsfs}
\usepackage{arydshln}
\usepackage{url,cite}
\usepackage{enumerate}
\usepackage{color}
\usepackage{pifont}

\newcommand{\blue}[1]{{\color{blue} #1}}

 \definecolor{mygray}{gray}{0.2}

\usepackage[normalem]{ulem}
\usepackage{soul}
\graphicspath{{images/}}

\usepackage{setspace}
\def\cals{\mathcal{S}}
\def\call{\mathcal{L}}
\def\calu{\mathcal{U}}

\def\caln{\mathcal{N}}

\def\C{C}
\def\M{\cals}
\def\R{\mathbb{R}}
\def\N{\mathbb{N}}
\def\bigdot{\boldsymbol{\cdot}}
\def\eref#1{{\rm (\ref{#1})}}
\def\n{{n}}
\def\t{{\beta}}

\def\to{\rightarrow}
\def\cp{{\mathrm{cp}}}
\def\H{{H}_{\cp}}
\DeclareMathOperator*{\arginf}{arg\,inf}
\DeclareMathOperator*{\esssup}{ess\,sup}

\def\Span{\mbox{span}}
\def\supp{\mbox{supp}}
\newtheorem{assumption}{Assumption}

\def\pn{\par\noindent}
\def\pf{\pn\textbf{Proof.}\enspace\ignorespaces}%
\def\qed{~\relax\ifmmode\hskip2em \Box
 \else\unskip\nobreak\hskip1em \hfill$\Box$
 \fi \newline}

\def\TM{\mathcal{T}_y{\M}}

\def\dS{{d_\M}}

\title{A kernel-based least-squares collocation method for surface diffusion}

\author{Meng Chen\footnotemark[1]$~^,$\footnotemark[2]
 	\and Ka Chun Cheung\footnotemark[3]$~^,$\footnotemark[4]
\and  Leevan Ling\footnotemark[4]
}

\footnotetext[1]{School of Mathematics and Computer Sciences, Nanchang University, Nanchang, China.}
\footnotetext[2]{Institute of Mathematics and Interdisciplinary Sciences, Nanchang University, Nanchang, China}
\footnotetext[3]{NVIDIA AI Technology Center (NVAITC), NVIDIA, USA.}
\footnotetext[4]{Department of Mathematics, Hong Kong Baptist University, Kowloon Tong, Hong Kong.}
\begin{document}

\maketitle
\begin{abstract}
There are plenty of applications and analysis for time-independent elliptic partial differential equations in the literature hinting at the benefits of overtesting by using more collocation conditions than the number of basis functions. Overtesting not only reduces the problem size, but is also known to be necessary for stability and convergence of widely used unsymmetric Kansa-type strong-form collocation methods.
We consider kernel-based meshfree methods, which is a method of lines with {\em collocation} and {\em overtesting} spatially, for solving parabolic partial differential equations on surfaces without parametrization.
In this paper, we extend the time-independent convergence theories for overtesting techniques to the parabolic equations on smooth and closed surfaces.
\end{abstract}

\section{Introduction}
For closed manifold $\M\subset\R^d$ and some given square integrable functions $f:\M\times[0, {T}]\to \R$ and $g:\M\to \R$,
we consider parabolic PDEs for some time-dependent surface scalar function $u:\M\times[0,T]\to\R$  in the form of
\begin{subequations}\label{eqPDE}
    \begin{align}
    \displaystyle
      \dot u(y,t) + \call_\M u(y,t) = f(y,t) &\quad \mbox{for }(y,t)\in\M\times[0,T] \label{eqPDE-a}
      \\
      u(y,0) = g(y)~~\, &\quad \mbox{for } y\in\M, \label{eqPDE-b}
    \end{align}
\end{subequations}
with a second-order uniformly elliptic operator in divergence form
\begin{equation}\label{equLs}
    \call_{\M } u(y,t):= - \nabla_{\M }\bigdot \big( A(y,t) \nabla_{\M } u(y,t) \big) \quad \mbox{for }(y,t)\in\M\times[0,T].
\end{equation}
Assumptions and definitions required to make sense of \eref{eqPDE} will be provide in Sect.~\ref{sec:Notations and preliminaries}.

We focus on strong-form collocation kernel-based spatial discretization. In \cite{Cheung+Ling-Kernembemethconv:18,Chen+Ling-Extrmeshcollmeth:20}, PDE \eref{eqPDE} were solved by the method of Rothe, in which the PDE is discretized first in time, and then by kernel-based collocation method in space.
In this paper, we focus on the method of lines that is theoretically more completed for parabolic PDEs on spheres \cite{Kunemund+NarcowichETAL-highmeshGalemeth:19} using Galerkin formulations and also in bulk domains, see \cite{Hon+SchabackETAL-meshKernmethline:14,Hon+Schaback-Diremeshkerntech:15}.
An advantage of using collocation method is that we can completely remove surface integrations from our algorithms (but not theories). Readers can find our theoretical study on the convergence of semi-discretized solution in Sect.~\ref{sec:Semi-discretized trial solution and convergence} (see Thm.~\ref{thm:conv}).
Algorithms for computing the fully discretized solution and the corresponding error analysis (see Thm~\ref{thm:conv 2}) were provided  in Sect.~\ref{sec:Fully discretized solution}. Numerical examples in Sect.~\ref{sec:Numerical examples} verify the convergence behaviour of the proposed methods and robustness in simulating solutions to Allen-Cahn equations.

\section{Notations and preliminaries}
\label{sec:Notations and preliminaries}
Let $\M\subset\R^d$ be a closed, connected, orientable, and complete Riemannian manifold of dimension $\text{dim}\,\M := \dS= d-1$.
We further suppose that  $\M$ is of class $\mathcal{C}^{m+1}$ for some  integer $m$
with bounded geometry and boundary regularity as required in the theories in \cite{Maerz+Macdonald-CALCSURFWITHGENE:12,Hangelbroek+NarcowichETAL-DireInveResuBoun:17}.
Under these assumptions, there exists $\delta>0$ such that Euclidean closest point retraction map
\begin{equation}\label{def:cp}
  \cp(x) := \arginf_{y\in\M} \|y-x\|_{ {\ell}^2(\R^d)}:\Omega\to \M
\end{equation}
is well-defined and $\mathcal{C}^{m}$-smooth in the narrow band domain
\begin{equation}\label{def:Omega}
  \Omega:= \{x \in\R^d\,:\,\|x-y\|_{ {\ell}^2(\R^d)}< \delta,\; y\in\M\}.
\end{equation}
Here, we use $\|x\|_{\ell^2(\R^d)}$ to denote the standard Euclidean 2-norm for any vector $x\in \R^d$.

\subsection{Differential operators  without parametrization}
\label{sec:Differential operators  without parametrization}
We adopt a set of parametrization-free definitions in \cite{Ruuth+Merriman-SimpEmbeMethSolv:08} for differential operators on $\M$.
Let $\n :\M\to \R^d$  be a smooth (column) vector field that spans the normal space.
We define the orthogonal projection matrix
\begin{equation}\label{def:P(y)}
  P(y) := I_d - n(y) n(y)^T,
\end{equation}
that projects onto the tangent space $\TM$ of $\M$ at $y\in\M$.
For any continuously differentiable scalar surface function $v:\M\to\R$, the surface gradient operator $\nabla_\M$ is defined by
\begin{equation}\label{def:grad}
  \nabla_\M v(y) := P(y)\nabla\big(v\circ\cp\big)(y)  = \nabla\big(v\circ\cp\big)(y),\quad y\in\M.
\end{equation}
For any continuously differentiable surface vector field $g:\M\to\R^{d}$, the surface divergence operator $\nabla_\M\bigdot$ is defined by
\begin{equation}\label{def:div}
  \nabla_\M \bigdot g(y) := \big(P(y) \nabla\big)\bigdot\big(g\circ\cp\big)(y) =  \nabla\bigdot\big(g\circ\cp\big)(y),\quad y\in\M.
\end{equation}
Note that these definitions are equivalent  \cite{Maerz+Macdonald-CALCSURFWITHGENE:12} to their intrinsic counterparts defined by local parametrization and Riemannian metric tensor.

\subsection{Diffusion tensor}
\label{sec:Diffusion tensor}
The diffusion tensor $A(y,t)$ in \eref{equLs} needs to satisfy some symmetric positive definiteness assumption when restricted to the tangent space $\TM$ of $\M$.

\begin{assumption}\label{assumption A1}
The symmetric surface intrinsic diffusion tensor $A:\M\times[0,{T}]\to \R^{d\times d}$ in \eref{equLs} satisfies,
for all $y\in\M$ and $t\in[0,T]$, that
\begin{subequations}\label{A1}
(1) $\TM$ is an invariant subspace of $A$, i.e.,
\begin{equation}\label{A1-a}
  \xi \in \TM \Rightarrow A(y,t)\xi \in \TM,
          \qquad \forall  \,\xi\in\TM,
\end{equation}
(2) $A(y,t)$ is uniformly positive definite on $\TM$, i.e., there exists $\nu>0$ such that
\begin{equation}\label{A1-b}
   \nu \| \xi\|_2^2 \leq  \xi^T A(y,t) \xi \leq \nu^{-1} \|\xi\|_2^2
         \qquad \forall  \,\xi\in\TM,
\end{equation}
and, (3) all components $a_{ij}(y,t):=[A(y,t)]_{ij}$, $1\leq i,j \leq d$, of $A$ are sufficient smooth with bounded derivatives.
\end{subequations}
\end{assumption}

Equivalently, \eref{A1-a}--\eref{A1-b} ensure the existence of a set of orthonormal eigenvectors   $\{\t_{i}(y,t)\}_{i=1}^{\dS}\subset \R^d$  of the diffusion tensor $A(y,t)$  that spans $\TM$
and all the associated eigenvalues are strictly positive.
Moreover,  the surface vector field $A(y,t)\nabla_\M u(y,t)$ is tangent to the surface and lies in $\TM$. Also, by \eref{A1-b} and  \cite[Thm. 4.2--4.3]{Maerz+Macdonald-CALCSURFWITHGENE:12}, we can recast the diffusion operator in \eref{equLs} by the Cartesian gradient and divergence operators  as
\begin{equation}\label{equLs-2}
	 \call_{\M }u(y,t)= - \nabla \bigdot \big( A(\cp(y),t) \nabla( u(\cp(y),t) ) \big) \quad \mbox{for }(y,t)\in\M\times[0,T],
\end{equation}
in terms of the $\cp$-mapping \eref{def:cp}. {Putting $A=I_d$ in \eref{equLs-2} yields the Laplace-Beltrami operator
	\begin{equation}\label{def:lap}
	\Delta_\M u(y,t) := \nabla_\M \bigdot (\nabla_\M u(y,t))=\Delta\big(u(\cp(y),t)\big),\quad  \mbox{for }(y,t)\in\M\times[0,T],
	\end{equation}
defined by $\cp$-mapping via \eref{def:grad} and \eref{def:div}.}

\subsection{{Derivatives and norms in Hilbert spaces}}
\label{sec:Sobolev norms}
Hilbert spaces  on complete manifolds with bounded geometry $\M$  are defined in \cite{Strichartz-AnalLaplcompRiem:83}  to be $H^k(\M):=(I-\Delta_\M)^{-k/2}L^2(\M)$ with $\Delta_\M$ in some
 equivalent  form to our $\cp$-version in \eref{def:lap}. This $H^k(\M)$ is norm equivalent \cite{Yoshida-SobospacRiemmani:92} to the Sobolev spaces characterized by localization  by an atlas
 {with open cover $\{U_i\}_i$ of $\M$,}
 and a subordinate  {$C_0^k(\M)$}  partition of unity  \cite[Def. 4.4]{Wloka-PartDiffEqua:87}.
We decompose  {$v\in W_2^k(\M)$,} whose precise definition is given in \cite[Thm. 2.13]{Wloka-PartDiffEqua:87},  as $v=\sum_i v_i$ with $\supp(v_i)\subset U_i$ and express $\M\cap U_i$ by, say,
$x_d = a_i(x_1,\ldots,x_{\dS})$  without loss of generality.
Let the multi-indexed  derivative with respect to the Cartesian coordinate $x_1,\ldots,x_{d}$   be
\begin{equation}\label{def:D^a}
  D_{d}^\alpha := \frac{\partial^{|\alpha|}}{\partial x_{1}^{\alpha_{1}}\cdots\partial x_{d}^{\alpha_{d}}},
  \quad |\alpha|={\alpha_{1}+\cdots+\alpha_{d}}
\end{equation}
and, we also define  $D_{\dS}^\alpha$  as in \eref{def:D^a}, but,  with respect to  $x_1,\ldots,x_{d-1}=x_{\dS}$.
Let the surface element with respect to variables $x_1,\ldots,x_{\dS}$ be
\[ 
  d\sigma_i = \bigg(
  1+\sum_{j=1}^{\dS} \Big( \frac{\partial a_i}{\partial x_j} \Big)^2 \bigg)^{\frac12}dx_1\cdots dx_{\dS}.
\] 
We define the non-standard $\H^k(\M)$-norm via $\cp$-operator by
\begin{equation}\label{def:W2k-norm}
\|v\|_{\H^k(\M)}^2 = \sum_{|\alpha|\leq k} \int_{\M} \big|
{D_\dS^\alpha} v \big|^2\,d\sigma
:= \sum_i \|v_i\|_{\H^k(U_i)}^2
\end{equation}
with
\begin{equation}\label{def:Ui-norm}
  \|v_i\|_{\H^k(U_i)}^2 :=
\sum_{|\alpha|\leq k} \int_{U_i} \Big|
\Big( D_{d}^\alpha (v_i\circ\cp)\Big)_{\big|x_d = a_i(x_1,\ldots,x_{\dS})} \Big|^2\,d\sigma_i,
  \quad k\in\N.
\end{equation}
We explicitly include the $\cp$-operator in definition to go along with \eref{equLs-2}.
For functions $f\in L^2(0,T;\H^k(\M))$, we adopt the following  space-time norms notation
\begin{equation}\label{L2Hknorm}
    \|f\|_{L^2(0,T;\H^k(\M))}
    =\big\|\,\| f \|_{H^k_{\cp}(\M)}\,   \big\|_{L^2(0,T)}
    =\Big( \int_{0}^{T}  \| f(\cdot,\tau) \|_{\H^k(\M)}^2 \,{d\tau} \Big)^{1/2}.
\end{equation}

\begin{lemma}\label{lem:norm equiv}
  For any $\C^{m+1}$ surface $\M$ satisfying the assumptions in Sect.~\ref{sec:Notations and preliminaries} and  all $k\leq m$,
  the $\H^k(\M)$-norm in \eref{def:W2k-norm} is equivalent to the $W_2^k(\M)$-norm defined  via atlas and subordinate partition of unity.
\end{lemma}

\pf
For $v \in W_2^k(\M)$, \cite[Sect. 4.2]{Wloka-PartDiffEqua:87}  shows that the $W_2^k(\M)$-norm (i.e., defined by an atlas and a subordinate partition of unity) is equivalent to a surface norm defined by means of coordinate invariant surface integrals as in \eref{def:W2k-norm}, but with localized counterparts
\begin{equation}\label{def:Ui-norm mod}
  \|v_i\|_{W_2^k(U_i)}^2 =
\sum_{|\alpha|\leq k} \int_{U_i} \Big|
D_{\dS}^\alpha \Big({v_i}_{\big|x_d = a_i(x_1,\ldots,x_{\dS})}\Big) \Big|^2\,d\sigma_i.
\end{equation}
Note that this $W_2^k(U_i)$-norm in \eref{def:Ui-norm mod} 
is exactly the one that can be deduced by using metric tensors of Riemannian manifold with the parameterized equation
\[
    r(x_1, x_2, \ldots, x_{\dS}) = \big(x_1, x_2, \ldots, x_{\dS},a_i(x_1,\ldots,x_{\dS}) \big)^T : \R^{\dS}\to\M\subset\R^d,
\]
which defines a set of basis of the tangent space
\[
    \{\t_j\}_{1\leq j\leq \dS}
    := \left\{ \frac{\partial r}{\partial{x_j}} \right\}_{1\leq j\leq \dS}
    = \left\{ \left[
      \begin{array}{c}
        {e}_j \\
        \frac{\partial a_i}{\partial{x_j} } \\
      \end{array}
    \right] \right\}_{1\leq j\leq \dS},
\]
where ${e}_j\in\R^{\dS}$ is the  $j$-th standard unit vector.
Without loss of generality, we assume ${\partial a_i}/{\partial x_j}\neq0$ for some $j$, or else $x_d = \text{const}$. In this case, the norms in \eref{def:Ui-norm} and \eref{def:Ui-norm mod} are trivially  equivalent.

Since $v_i = v_i\circ\cp$ for all  $y=(x_1,\ldots,x_d)\in\M$,
we now work on the function in the integrand of \eref{def:Ui-norm mod} and define
$
    V(x_1,\ldots,x_{\dS}) :=\big(v_i\circ\cp\big) 
    _{\big| x_d = a_i(x_1,\ldots,x_{\dS})}:\R^{\dS}\to\R.
$
By  implicit differentiation,
we know that the first derivatives were connected by
\[
    \frac{\partial}{\partial x_j} V
     = \frac{\partial}{\partial x_j}( v_i \circ \cp)
     + \frac{\partial  a_i}{\partial x_j}\frac{\partial}{\partial x_d}( v_i \circ \cp),
     \qquad j=1,\dots,\dS.
\]
Rewriting in matrix form yields
\begin{eqnarray}
    \nabla_{{\dS}} V(x_1,\ldots,x_{\dS})
    &=& \big[ \t_1,\dots,\t_{\dS} \big]^T \nabla_{ {d}} ( v_i \circ \cp)(y).
    \label{eq:gradV}
\end{eqnarray}

{The norm equivalency follows immediately from the identity \eref{def:Ui-norm mod}  by the chain rule, and by the fact that $\nabla_{{d}} (v_i\circ\cp)(y)\in\TM$, we can uniquely express
\begin{equation}\label{eq:gradV2}
    \nabla_{{d}} (v_i\circ\cp)(y) = \sum_{j=1}^{\dS}\gamma_j\, \t_j
\end{equation}
with some coefficients $\gamma=[\gamma_1, \ldots, \gamma_{\dS}]^T\in \R^{\dS}$.
\qed}

\section{Semi-discretized trial solution and its convergence}
\label{sec:Semi-discretized trial solution and convergence}
We can apply the method of lines to discretize the surface diffusion equation \eref{eqPDE} spatially by some kernel-based trial space.
In this paper, we focus on manifold kernels
that can be obtained by restricting \cite{Fuselier+Wright-ScatDataInteEmbe:12,narcowich2007approximation} some global, symmetric positive definite, and Sobolev space $H^{m+1/2}(\R^d)$ reproducing \cite{Matern} kernels ${\Phi}_{m+1/2}: \R^d \times \R^d \rightarrow \R$ to $\M$.
Fourier transforms of such kernels $\hat{\Phi}_\tau$  decay like
\begin{equation}\label{kernelFour}
c_1(1+\|{\omega}\|_2^2)^{-{(m+1/2)}}\le  \hat{\Phi}_{m+1/2}({\omega})
\leq c_2(1+\|{\omega}\|_2^2)^{-{(m+1/2)}} \quad \mbox{for all }{\omega}\in\R^d,
\end{equation}
for some constants $0<c_1\leq c_2$.
Define the manifold kernels $\Psi_m:\M\times\M\to\R$ by
\begin{equation}\label{eq:MKernels}
  \Psi_m(\cdot,\cdot) := {\Phi}_{m+1/2}(\cdot,\cdot)_{|\M\times\M}.
\end{equation}
Then, for any $m>\dS/2=(d-1)/2$, this manifold kernel $\Psi_m$ reproduces $H^{m}(\M )$, see  \cite{Fuselier+Wright-ScatDataInteEmbe:12}. 
In practice, one can use the standard Whittle-Mat\'{e}rn-Sobolev kernels \cite{Matern} or the  Wendland compactly supported kernels \cite{Wendlandfun} with smoothness order $m+1/2$, as is, for implementation.

\subsection{Trial spaces}
Let  $Z=\{z_1, \ldots, z_{n_Z}\} \subset \M$ be the set of trial centers. The {fill distance} $h_{Z}$ and the {separation distance} $q_{Z}$ are defined respectively as
\begin{equation}\label{def_h_q}
h_{Z}:=\sup_{\zeta\in \M} \inf_{\eta\in Z} \text{dist}({\zeta},{\eta}) \quad \mbox{and}\quad q_{Z}:=\frac{1}{2}\inf_{{{z}_i\neq {z}_j \in Z}}\text{dist}({z}_i,{z}_j).
\end{equation}
The {mesh ratio} is $\rho_{Z} =h_Z/q_Z \geq 1$. We say that surface points in $Z$ are quasi-uniform if, as $n_Z:=|Z|$ increases, $Z$  satisfies
$q_{Z}\leq h_{Z}= \rho_Z q_{Z} \leq \rho q_Z,$
for some constant $\rho>0$ independent of $n_Z$.
Because of equivalence between them on  smooth and compact surfaces  \cite[Thm. 6]{Fuselier+Wright-ScatDataInteEmbe:12},
we can simply use Euclidean {distance} as $\text{dist}(\cdot,\cdot)$ instead of geodesic distance.

We work on finite-dimensional trial spaces $\calu_Z$,  that is a span of the translation-invariant $\Psi_m$ in \eref{eq:MKernels} to $Z$,  defined by
\begin{equation}\label{Uzfinite}
\calu_{Z,\M ,\Psi_m}:= \Span\{\Psi_m\left(\cdot\, ,z_j\right)\;|\;z_j\in Z\}.
\end{equation}
Any (time-dependent) trial function $u\in \calu_{Z,\M ,\Psi_m}$ can be expressed by a linear combination with a set of unknown coefficients (function of time) $\{\lambda_j\}$
\begin{equation}\label{eq:trial sp0}
u_Z  = \sum_{z_j\in Z}\lambda_j  \Psi_m(\cdot,z_j)
=: \Psi_m( \cdot, Z) \lambda_Z .
\end{equation}
The associated reproducing kernel Hilbert space (a.k.a native space) $\mathcal{N}_{\Psi_m}(\M)$-norm is given by
\begin{equation}\label{eq:native_norm}
\|u_Z\|_{\mathcal{N}_{\Psi_m}(\M)}^2
=\sum_{z_i,z_j\in {Z}} \Psi_m(z_i,z_j) \lambda_i \lambda_j
=: \lambda_Z^T\Psi_m(Z, Z)\lambda_Z.
\end{equation}
By \cite{Buhmann-Radibasifunc:03,Wendland-Erroestiintecomp:98} and Lem.~\ref{lem:norm equiv}, we know that $\|u_Z\|_{\mathcal{N}_{\Psi_m}(\M)}\sim\|u_Z\|_{W_2^m(\M)}\sim\|u_Z\|_{\H^m(\M)}$ as defined in \eref{def:W2k-norm}. Throughout the paper, we use  $\H^m(\M)$-norm in place of native space norm for simplicity and norm equivalency will be taken care by generic constants.

\subsection{Strong form collocation and overtesting}
Let $Y = \{y_1,\ldots, y_{n_Y}\} \subset \M$,  $n_Y>n_Z$,   be a sufficiently dense set of quasi-uniform collocation points.
We use strong form collocation at $Y$ as \emph{test conditions} in order to identify a numerical solution form the trial space $\calu_{Z,\M ,\Psi_m}$.
Let  the discrete $\ell^2(Y)$-norm be
\begin{equation}\label{eq:ell2Y}
 \| f \|_{\ell^2(Y)}^2  := \sum_{y_i\in Y} |f(y_i)|^2.
\end{equation}
For quasi-uniform $Y$, we have $
	\left(
\displaystyle	\int_\M f(y)dy
	\right)^2
= \calo\big( h_Y^{\dS} \| f \|_{\ell^2(Y)}^2 \big)
$. For some regularization parameter $\alpha$, we now formally define the {\emph{time-dependent semi-discretized trial solution} to \eref{eqPDE} in the time-dependent trial space   $L^2(0,T;\calu_{Z,\M ,\Psi_m})$, whose norm is defined similarly to \eref{L2Hknorm} but with finite dimensional trial spaces \eref{Uzfinite} in place of $\H^k(\M)$, by following optimization problem}
 \begin{equation}\label{eq:un}
      \begin{array}{l}
    { u_{Z,\alpha} }:= \displaystyle
{h_Y^{\dS-2}}\!\!\!\!
\arginf_{
  u,\dot u \in L^2(0,T;\calu_{Z,\M ,\Psi_m})
}
\Big({h_Y^2}\int_{0}^{{T}}\|\dot{u}{(\cdot,\tau)}+\call_\M u{(\cdot,\tau)}-f {(\cdot,\tau)}\|_{\ell^2(Y)}^2 d{\tau}
      \\\qquad\qquad\qquad\qquad\qquad\qquad\qquad
      + \|u(\cdot,0)-g{(\cdot)}\|_{\ell^2(Y)}^2+ \alpha^2\|u(\cdot,0)\|^2_{\H^m(\M)}\Big),
      \end{array}
\end{equation}
for some given initial condition  {to be discussed later in Sect.~\ref{sec:Regularized initial condition}  rigorously}. Note that the scaling factor $h_Y^{\dS-2}$ in \eref{eq:un} is for the sake of convergence analysis in Thm. \ref{thm:conv} below and it is not necessary in computations.

Putting the trial function in the form of \eref{eq:trial sp0} into the PDE residual, $\dot{u}+\call_\M u-f$, and evaluating at $y_i\in Y$ yields the following  {expressions}
\[
    \sum_{z_j\in Z}\Psi_m(y_i,z_j) \frac{d}{dt}\lambda_j(t)
    +\sum_{z_j\in Z}[\call_\M\Psi_m](y_i,z_j) \lambda_j(t)
    -f(y_i,t)
    \quad \text{for } t\in(0,T],
\]
or in an overdetermined matrix as
\[
    \Psi_m(Y,Z)\dot{\lambda}_Z(t) + [\call_\M\Psi_m](Y,Z)\lambda_Z(t) - f(Y,t)
    \quad \text{for } t\in(0,T].
\]
 {Since each (time dependent $L^2$-integrable) trial function in $\calu_{Z,\M ,\Psi_m}$ can be uniquely identified by a vector of unknown (time-dependent) coefficient  (in $(L^\infty[0,{T}])^{|{Z}|}$) in $\R^{|Z|}$,}   we can recast
\eref{eq:un} in terms of the unknown coefficient functions as %
 \begin{equation}\label{def:uN-a}
      \begin{array}{l}\displaystyle
   \lambda_{Z,\alpha} {(\cdot)} := \arginf_{
   \lambda,\dot\lambda\in (L^2\cap L^\infty[0,{T}])^{|{Z}|}
   }
   \int_{0}^{{T}}
  \big\| \Psi_m(Y,Z)\dot\lambda{(\tau)}
  \\
   \qquad \qquad \qquad \qquad \qquad \qquad + [\call_\M \Psi_m](Y,Z)\lambda{(\tau)}-f(Y,{\tau}) \big \|^2_{\ell^2(\R^{|Y|})}    d{\tau},
      \end{array}
\end{equation}
subject to  {the same initial conditions accompanying \eref{eq:un} in terms of  unknown coefficients}.

In the  {remaining} of this section, we will prove convergence estimates by proving a regularity estimate for surface diffusion \eref{eqPDE-a} in Sect.~\ref{sec:Regularity estimate},
an error estimate of a regularized initial condition $g_{Z,\alpha}$ in Sect.~\ref{sec:Regularized initial condition},
 and,
finally in Sect.~\ref{sec:Convergence estimate},
the convergence of the trial function $u_\alpha(\cdot,t):=\Psi_m(\cdot,Z) \lambda_{Z,\alpha}(t)$, which is defined via the solution of \eref{def:uN-a}, to the exact solution $u^*(\cdot,t)$.
We summarize the main result in Thm.~\ref{thm:conv}.

\subsection{Regularity estimates  for surface diffusion}
\label{sec:Regularity estimate}
We will walk through the details in proof to identify all necessary assumptions so that
\begin{eqnarray}
  &&\esssup_{0\leq {\tau}\leq T}\|{u{(\cdot,\tau)}}\|_{\H^1(\M)}+\|{u}\|_{L^2(0,T;\H^{\mu}(\M))}
  +\|{\dot u}\|_{L^2(0,T;\H^0(\M))}
\nonumber\\
  &&\qquad\qquad\qquad\qquad\qquad
  \leq C \Big(\|f\|_{L^2(0,T;\H^0(\M))}  + \|g\|_{\H^1(\M)}\Big),
    \mbox{\quad for $\mu=1,2$},\label{eq:energy est}
\end{eqnarray}
holds for the solution ${u}$ to the surface diffusion in \eref{eqPDE} with Sobolev norms defined in Sect.~\ref{sec:Sobolev norms}.
After we prepare the manifold versions of all required (in)equalities, regularity estimate \eref{eq:energy est} can be shown by applying some standard arguments.  We begin with the Green's first identity.

\begin{lemma}\label{lem:Green}
   Let $\M\subset\R^d$ be a $\C^3$-smooth closed manifold  and  $v,w:\M\to\R$   some $\C^2$-smooth scalar value functions.
   Then we have
   \[
    \int_\M v\Delta_\M (A w) \,d\sigma= -\int_\M  \nabla_\M v \bigdot(A\nabla_\M w)\,d\sigma.
    \]
    for any diffusion tensor $A$ satisfying Assumption~\ref{assumption A1}.
\end{lemma}
\pf Without loss of generality, we can ignore the time-dependency in the diffusion tensor $A$ in this proof.  We $\cp$-extend all functions to the narrow band domain $\Omega\supset \M$ in \eref{def:Omega} and denote
$v_\cp:=v\circ\cp,$ $ w_\cp:=w\circ\cp$,
    and  $ A_\cp:=A\circ\cp.$
We start with a vector calculus inequality\footnote{$\nabla\bigdot(f\boldsymbol{F}) = f\nabla\bigdot \boldsymbol{F} + (\nabla f)\bigdot \boldsymbol{F}$}:
\begin{equation}\label{eq:subgreen}
  \nabla\bigdot( v_\cp A_{cp}\nabla w_\cp)
= v_\cp \nabla\bigdot( A_{cp}\nabla w_\cp)
+ \nabla v_\cp\bigdot (A_{cp}\nabla w_\cp ).
\end{equation}
All smoothness requirements together ensure \eref{eq:subgreen} is well-defined and $v_\cp A_{cp}\nabla w_\cp$ is a smooth vector field.
Applying divergence theorem {for every smooth vector field $X$} on a closed manifold \cite[Thm. 1.1]{Pigola+Setti-Globdivetheononl:14}, i.e., $\int_\M X \,d\sigma=0$, to \eref{eq:subgreen} yield
\begin{eqnarray*}
    0 &=& \int_\M \nabla\bigdot( v_\cp A_{cp}\nabla w_\cp) \,d\sigma
= \int_\M v_\cp \nabla\bigdot( A_{cp}\nabla w_\cp)
+ \nabla v_\cp\bigdot (A_{cp}\nabla w_\cp ) \,d\sigma.
\end{eqnarray*}
The proof is completed after simplification by definitions \eref{def:grad} and \eref{def:div}. \qed

Using Lem.~\ref{lem:Green}, we  test \eref{eqPDE-a} by $u$ to obtain
\[
    \frac12\frac{d}{dt}\int_\M u^2(\cdot, t) \,d\sigma +  \int_\M \|(\nabla_\M u)(\cdot, t)\|_A^2 \,d\sigma = \int_\M (uf)(\cdot, t) \,d\sigma,
\]
where $\|\nabla_\M u\|_A^2 :={ \nabla_\M u \bigdot(A\nabla_\M u)}$.
Following  some standard arguments, i.e., integrate over time
\[
\frac12 \left(\int_\M u^2 \,d\sigma - \int_\M g^2\,d\sigma\right) +  \int_0^t \int_\M \|\nabla_\M u\|_A^2 \,d\sigma d\tau = \int_0^t \int_\M uf \,d\sigma d\tau,
\]
and apply Cauchy inequality $ab\leq \frac{\epsilon}{2}a^2+\frac{1}{2\epsilon}b^2$ with $\varepsilon = {2T}$, yields
\begin{eqnarray}
\int_0^t \int_\M uf \,d\sigma d\tau
&\leq&
{T} \int_0^T \int_\M f^2 \,d\sigma d\tau
+\frac{1}{4} \esssup_{0\leq \tau\leq T} \int_\M u^2{(\cdot,\tau)} \,d\sigma\label{eq:CI}.
\end{eqnarray}
We obtain an energy estimate
\begin{equation}\label{eq:EnergyEst1}
	 \esssup_{0\leq \tau\leq T}  \int_\M u^2{(\cdot,\tau)} \,d\sigma + \int_0^T \int_\M \|\nabla_\M u\|_A^2 \,\,d\sigma d\tau
	 \leq C\Big( \int_\M g^2\,\,d\sigma  +  \int_0^T \int_\M f^2 \,\,d\sigma d\tau	\Big),
\end{equation}
for some constant $C>0$ independent of the solution $u$ to \eref{eqPDE}. Moreover, due to {\eref{A1-b}  in} Assumption~\ref{assumption A1} , we know that
{
\[
\int_0^T\int_\M \|\nabla_\M u\|_A^2 \,d\sigma d\tau =\int_0^T \int_\M (\nabla_\M u)^T A(\nabla_\M u)\,d\sigma d\tau
\]
is equivalent to $\|\nabla_\M u\|^2_{L^2(0,T;\H^0(\M))}$} with constants depending on the eigenvalues of $A$ restricted to the tangent space $\TM$.

To get an {improved} regularity, we follow \cite[Thm.5, Ch.7]{Evans-Partdiffequa:98} and suppose that the diffusion tensor $A$ does not depend on $t$. Testing the PDE by $\dot{u}$ yields
\[
     \int_\M \dot{u}^2\,d\sigma  +\frac{d}{dt} \int_\M \frac12\|\nabla_\M u\|_A^2\,d\sigma = \int_\M \dot{u} f \,d\sigma,
\]
which, after an integration in time, will lead us to the estimate
\begin{equation}\label{eq:EnergyEst2}
 \esssup_{0\leq \tau\leq T} \|\nabla_\M u{(\cdot,\tau)}\|_A^2 + \int_0^T  \int_\M \dot{u}^2\,d\sigma  d\tau
	 \leq C\Big( \int_\M \|\nabla_\M g\|_A^2\,\,d\sigma  +  \int_0^T \int_\M f^2 \,\,d\sigma d\tau	\Big).
\end{equation}
For strong solution, the last required estimate comes form the fact that $\call_\M u(\cdot,\tau)\in L^2(\M)$ for a.e. $0\leq \tau\leq T$ and, by elliptic regularity \cite[Lem. 2.1 \& 3.1]{Chen+Ling-Extrmeshcollmeth:20} of the surface operator $\call_\M$ and \eref{eq:EnergyEst2}, we have
\begin{eqnarray}
     \int_0^T   \| u\|^2_{\H^2(\M)} \,d\tau
    &\leq&
    C   \int_0^T  \| \call_\M u \|^2_{\H^0(\M)}  \,  d\tau
    \nonumber  \\
    &\leq &
    C  \int_0^T   \Big( \| f\|^2_{\H^0(\M)}+ {\| \dot{u} \|^2_{\H^0(\M)} } \Big) \, d\tau
    \nonumber  \\
    &\leq &  C  \int_0^T    \| f\|^2_{\H^0(\M)} \, d\tau + {\int_\M \|\nabla_\M g\|_A^2\,\,d\sigma}.
    \label{eq:EnergyEst3}
\end{eqnarray}
Using the trivial bound $\|\,\bigdot\,\|^2_{L^2(0,T;\H^0(\M))} \leq C_\M   \esssup_{0\leq \tau\leq T} \|\,\bigdot\,\|^2_{\H^0(\M)) }$ to combine \eref{eq:EnergyEst1}--\eref{eq:EnergyEst3} completes the proof of the regularity estimate in \eref{eq:energy est}.

\begin{lemma}\label{lem:reg est}
  Assume the data functions in the surface diffusion equation \eref{eqPDE} with time-independent diffusion tensor are sufficiently smooth,
  $f\in L^2(0,T;L^2(\M))$, and $g\in H^1(\M)$,
  then the regularity estimate \eref{eq:energy est} holds on the unique strong solution
  $u^*\in (L^2\cap L^\infty)(0,T;H^2(\M))$ and $\dot{u}^*\in L^2(0,T;L^2(\M))$ to \eref{eqPDE}.
\end{lemma}

\subsection{Regularized discrete least-squares initial condition}
\label{sec:Regularized initial condition}
	
{In this section, we focus on regularized least-squares initial condition.} We show that the solution to the following time independent regularized approximation problem
\begin{equation}\label{def:gzalpha}
    g_{Z,\alpha}:=  \arginf_{u\in \calu_{Z,\M ,\Psi_m}}
    \|u -g\|_{\ell^2(Y)}^2
    + \alpha^2\|u\|^2_{\H^m(\M)}
\end{equation}
in \eref{eq:un} with some appropriate $\alpha$ is $\H^1(\M)$-convergent so that, later in our convergence analysis, we can show consistency by some bounds on $\|g_{Z,\alpha}-g\|_{\H^1(\M)}$ in the right handed side of \eref{eq:energy est}.
To do so, we extend the convergence estimate in \cite{2019Discrete} for flat geometry  to manifolds.
{When dealing with initial condition $g$ in \eref{eqPDE-b}  that does not suffer from Runge \cite{Bengt2007The} and Gibbs \cite{fornberg2011gibbs} phenomenon in radial basis function interpolation, one can use the interpolant of $g$ from $\calu_{Z,\M ,\Psi_m}$  in \eref{Uzfinite} as initial condition for defining the semi-discretized solution in \eref{eq:un} and \eref{def:uN-a}. This modification will not affect the convergence analysis below.}

{
	Note that the initial condition in terms of coefficients corresponding to \eref{def:gzalpha}  could be  determined by
{\begin{equation}\label{def:uN-b}
	  \lambda_{Z,\alpha}(0):=
			\arginf_{\lambda\in \R^{|Z|}} {\Big\{}\big\|\Psi_m(Y,Z)\lambda-g(Y) \big\|_{\ell^2(\R^{|Y|})}^2+ \alpha^2 \lambda^T\Psi_m(Z,Z)\lambda{\Big\}}.
		\end{equation}}
For the second initial condition $\dot\lambda_{Z,\alpha}(0)$, one can differentiate  \eref{def:uN-a}	or by other means, say
	$\dot{u}_{Z,\alpha}(\cdot,0) = \Pi\big(  f(\cdot,0)-\call_\M g \big)$
	with some appropriate projection $\Pi$ to the trial space.
	Then, the unique solvability of \eref{eq:un} is guaranteed by theorems in calculus of variations.
}

\begin{lemma}\label{lem:gz}
  For any $m \geq \lfloor1+\dS/2\rfloor$ and $\alpha\geq 0$, suppose that the sets  $Y, \,Z\subset \M$ of sufficiently dense discrete data points satisfy \eref{eq:alpha*}. Let
  $g_{Z,\alpha}\in \calu_{Z,\M ,\Psi_m}$ be the regularized approximant in \eref{def:gzalpha} for $g\in\H^m(\M)$. Then,
  \[
     \|g_{Z,\alpha}-g\|_{\H^1(\M)}^2 \leq C \Big(
     h_Y^{-2} h_Z^{2m-\dS}
     +h_Y^{\dS-2}\alpha^2
     +h_Z^{2m-3}
     \Big) \|g\|_{\H^m(\M)}^2
  \]
  holds for some constant $C$ independent of $g$.
\end{lemma}
When $\alpha \geq \alpha^*:= h_Y^{m-\dS/2}$, we have
\[
     \|g_{Z,\alpha^*}-g\|_{\H^1(\M)}^2 \leq C \Big(
     h_Y^{-2} h_Z^{2m-\dS}
     +h_Y^{2m-3}
     +h_Z^{2m-3}
     \Big) \|g\|_{\H^m(\M)}^2,
\]
without any extra denseness restriction imposed on $Z$ and $Y$.
\pf
We begin with a sampling inequality for manifold functions $g\in\H^m(\M)$, see {the second inequality on boundary in} \cite[Lem. 3.1]{Cheung+LingETAL-leaskerncollmeth:18} 
(with $\mathrm{m}=m+1/2$ and $\mathrm{s}=3/2$); for any discrete set $Y\subset\M$ with sufficiently small mesh norm $h_Y$, the following holds 
\begin{equation}\label{eq:SampIneqG}
\|g\|_{\H^1(\M)}^2 
\leq C h_Y^{d_S-2} \Big( \|g\|_{ \ell^2(Y)}^2 + h_Y^{2m-d_S} \|g\|_{\H^{m}(\M)}^2\Big) \end{equation}

Let $I_Z g$ be the interpolant of $g$ from the trial space $\calu_{Z,\M ,\Psi_m}$.
We have the following convergence estimate in  \cite[Cor.13]{Fuselier+Wright-ScatDataInteEmbe:12}
and \cite[Sect.3]{Chen+Ling-Extrmeshcollmeth:20}:
\begin{equation}\label{eq:Izg conv}
	\|I_Z g - g\|_{W_q^k(\M)} \leq C  h_Z^{m-k-\dS(1/2-1/q)_+} \| g\|_{\caln_{\Psi_m}(\M)}
\end{equation}
for 
$0\leq k\leq \lceil m-\dS(1/2-1/q)_+\rceil -1$.
Using \eref{eq:Izg conv} with $k=1$, it is sufficient to demonstrate the convergence $g_{Z,\alpha}\to I_Zg$ in the  trial space $\calu_{Z,\M ,\Psi_m}$, within which
we know a Bernstein's inverse inequality \cite[Thm. 10]{Hangelbroek+NarcowichETAL-DireInveResuBoun:17}
\begin{equation}\label{eq:invineq}
 \| u\|_{\H^k(\M)} \leq C h_Z^{-k} \|u\|_{L^2(\M)}, \qquad 0\leq k \leq m
\end{equation}
holds for all $u\in\calu_{Z,\M ,\Psi_m}$.

For any $\alpha\geq0$,  
we use \eref{eq:SampIneqG} and \eref{eq:invineq} on $e:=g_{Z,\alpha}-I_Zg\in \calu_{Z,\M ,\Psi_m} \subset \H^m(\M)$ to establish a chain of upper bounds:
\begin{eqnarray}
\|e\|_{\H^1(\M)}^2
&\leq&
C h_Y^{\dS-2} \Big(
\|e\|_{ \ell^2(Y)}^2
+  h_Y^{2m-\dS} \|e\|_{\H^{m}(\M)}^2 \Big)
\nonumber
\\
&=&
C h_Y^{\dS-2} \Big(
\|e\|_{ \ell^2(Y)}^2
+  \alpha^2 \|e\|_{\H^{m}(\M)}^2 + (  h_Y^{2m-\dS}-\alpha^2)_+ \|e\|_{\H^{m}(\M)}^2 \Big)
\nonumber
\\
&\leq&
C h_Y^{\dS-2} \Big(
\|e\|_{ \ell^2(Y)}^2
+  \alpha^2 \|e\|_{\H^{m}(\M)}^2 + (  h_Y^{2m-\dS}-\alpha^2)_+ h_Z^{-2m}\|e\|_{L^2(\M)}^2 \Big)
\nonumber
\\
&\leq&
C h_Y^{\dS-2} \Big(
\|e\|_{ \ell^2(Y)}^2
+  \alpha^2 \|e\|_{\H^{m}(\M)}^2  \Big),
\label{H1sampling}
\end{eqnarray}
either if $\alpha \geq \alpha^*:= h_Y^{m-\dS/2}$, or
under the {denseness} and regularization constraint
\begin{equation}\label{eq:alpha*}
  C h_Y^{\dS-2}(  h_Y^{2m-\dS}-\alpha^2)_+ h_Z^{-2m}   <\frac12.
\end{equation}
{Note that similar denseness conditions of two point sets were also required in all previous works \cite{Cheung+LingETAL-leaskerncollmeth:18,Cheung+Ling-Kernembemethconv:18,Chen+Ling-Extrmeshcollmeth:20} for  stability estimates to hold in theory. Similar to all previous observations, this sufficient condition is not numerically necessary as we will soon see in the later numerical examples, where linear ratio of oversampling is numerically sufficient. There are rooms for further theoretical improvement and we leave this to future works.}

Add-in subtract-out the function $g$ in \eref{H1sampling};
a direct consequence of the optimality in \eref{def:gzalpha} is that
\begin{equation}
  \|e\|_{\H^1(\M)}^2 
  \leq 2C h_Y^{\dS-2} \Big(
\|I_Zg - g\|_{ \ell^2(Y)}^2
+  \alpha^2 \|I_Zg-g\|_{\H^{m}(\M)}^2  \Big).
\end{equation}
A scattered zero lemma
\cite[{Lem}.10]{Fuselier+Wright-ScatDataInteEmbe:12} suggests that
\begin{eqnarray*}
\|I_Zg-g\|_{ \ell^2(Y)}^2
&\leq& n_Y  \|I_Zg-g\|_{ L^\infty(\M) }^2
\leq  C  n_Y   h_Z^{2m-\dS} \|I_Zg-g\|_{\H^m(\M)}^2,
\end{eqnarray*}
provided $m\geq\lfloor1+\dS/2\rfloor$.
Because $I_Zg$ and $I_Zg-g$ are mutually orthogonal with respect to the native space norm, i.e., $\H^m(\M)$-norm, we know that both $\|I_Zg\|_{\H^m(\M)}$ and $\|I_Zg-g\|_{\H^m(\M)}\leq \|g\|_{\H^m(\M)}$.
Under the assumption that $Y$ is quasi-uniform and therefore $n_Y \leq c q_Y^{-\dS} \leq C h_Y^{-\dS}$, we arrive at
\begin{equation}\label{|IZg-g|_Y}
  \|I_Zg-g\|_{ \ell^2(Y)}^2 \leq C h_Y^{-\dS} h_Z^{2m-\dS} \|g\|_{\H^m(\M)}^2
\end{equation}
and complete the proof.
\qed

\subsection{Convergence estimate for semi-discretized solution}
\label{sec:Convergence estimate}~

Define a residual functional $\mathcal{E}_\mu:{H^1}(0,T;\H^m(\M))\to \R$ using the left {hand} side of the regularity estimate \eref{eq:energy est} for $\mu=1,2$ by
\begin{equation}\label{eq:errorfunctional}
  \mathcal{E}_\mu {[u]} :=
  \esssup_{0\leq \tau\leq T}\|u{(\cdot,\tau)}\|_{\H^1(\M)}^2
  +\|u\|_{L^2(0,T;\H^\mu(\M))}^2
  +\|\dot u\|_{L^2(0,T;\H^0(\M))}^2,
\end{equation}
and restate \eref{eq:energy est} as
\begin{equation}\label{eq:errorfunctional bound}
 \mathcal{E}_\mu {[u]}  \leq C \Big(
   \|\dot{u} + \call_\M u \|_{L^2(0,T;\H^0(\M))}^2
    + \|u(\cdot,0)\|_{\H^1(\M)}^2
 \Big) .
\end{equation}
We begin by showing that there is a good comparison function in the trial space.
For any $t\in[0,T]$, let $s(\cdot,t) := I_Z u^*(\cdot, t)$ be the interpolant  of the solution $u^*$ to the PDE~\eref{eqPDE} from the trial space $\calu_{Z,\M ,\Psi_m}$. Because
\begin{equation}\label{eq:s}
    s(z_j,t) = I_Z u^*(z_j, t) = \Psi_m(z_j,Z)[\Psi_m(Z,Z)]^{-1} u^*(Z, t)
    \quad\text{for all $z_j\in Z$ and $t\in[0,T]$},
\end{equation}
differentiating with respect to $t$ shows that
\[
   \dot s(z_j,t) =  \Psi_m(z_j,Z)[\Psi_m(Z,Z)]^{-1} \dot u^*(Z, t)
    \quad\text{for all $z_j\in Z$ and $t\in[0,T]$},
\]
and therefore $\dot{s}(\cdot,t) \in \calu_{Z,\M ,\Psi_m}$ is the unique interpolant of the first order time derivative $\dot{u}^*(\cdot,t)$.
Using \eref{eq:errorfunctional bound} and standard interpolation theories \cite[Cor. 13]{Fuselier+Wright-ScatDataInteEmbe:12},
we get the following error estimate:
\begin{eqnarray*}
    &&\mathcal{E}_\mu  {[s-u^*]}
   \nonumber  \\
   &&\quad \leq C \Big(
   \|\dot{s} -\dot{u}^*\|_{L^2(0,T;\H^0(\M))}^2
   + \| \call_\M (s-u^*) \|_{L^2(0,T;\H^0(\M))}^2
   + \|(s-u^*)(\cdot,0)\|_{\H^1(\M)}^2 \Big)
   \nonumber  \\
   &&\quad \leq C h_Z^{2m-4} \Big(
       \| \dot{u}^*\|_{L^2(0,T;\H^{m-2}(\M))}^2
   +   \| u^*\|_{L^2(0,T;\H^m(\M))}^2
   +   h_Z^{2}\| u^*(\cdot,0) \|_{\H^m(\M)}^2   \Big).
   \label{E(s-u*)}
\end{eqnarray*}

Our next task is to show that the numerical solution defined by \eref{eq:un} converges to the comparison function $s$, i.e., $\mathcal{E}_\mu {[u_{Z,\alpha} - s]}\to 0$ as $h_Z\to0$, in the trial space. Below are two essential stability estimates before we can study consistency. The first comes from \eref{H1sampling}.

\begin{corollary}\label{lem:stability2}
{Let $m \geq \lfloor  1+\dS/2\rfloor$}
and $\alpha\geq 0$.  If $\alpha < \alpha^*:= h_Y^{m-\dS/2}$, further suppose that the sets  $Y, \,Z\subset \M$ of sufficiently dense discrete data points satisfy \eref{eq:alpha*}. Then,
\[
\|u\|_{\H^1(\M)}^2 \leq
C h_Y^{\dS-2} \Big(
\|u\|_{ \ell^2(Y)}^2
+  \alpha^2 \|u\|_{\H^{m}(\M)}^2  \Big),
\]
holds for all trial function $u \in \calu_{Z,\M ,\Psi_m}$.
\end{corollary}

\begin{lemma}\label{lem:stability}
{Let $m\geq \lfloor3+\dS/2\rfloor$}.
Suppose that the sets  $Y, \,Z\subset \M$ of sufficiently dense discrete data points satisfy \eref{eq:desnseness stability}. Let $\mathcal{E}_\mu {[u]} $ be the residual functional defined in \eref{eq:errorfunctional}. Then, the estimate
\[
\mathcal{E}_\mu {[u]}
\leq  C \Big(  h_Y^{\dS}  \int_{0}^{T} \|\dot{u} {(\cdot, \tau)}+   \call_\M u{(\cdot, \tau)}\|_{ \ell^2(Y)}^2 \,d{\tau}
    + \|u(\cdot,0)\|_{\H^1(\M)}^2
 \Big)
\]
holds for all trial functions $u,\dot{u}\in L^2(0,T; \calu_{Z,\M ,\Psi_m} )$ and $\mu=1,2$.
\end{lemma}
\pf
For any $t\in[0,T]$, applying the sampling inequality in \cite[Lem. 3.1]{Cheung+LingETAL-leaskerncollmeth:18} (with $\mathrm{m}=m+1/2$ and $\mathrm{s}=1/2$) to $\dot{u}(\cdot, t)+\call_\M u(\cdot, t) \in \H^{m-2}(\M)$
for any time $t\in[0,T]$
in \eref{eq:errorfunctional bound} at the set $Y$, then inverse inequality \eref{eq:invineq} and the estimate $\| \call_\M u \|_{\H^{m-2}(\M)}\leq \|   u \|_{\H^{m}(\M)}$ due to boundedness of PDE coefficients,  yields
\begin{eqnarray}
   &&\| \dot{u}+\call_\M u \|_{\H^0(\M)}^2
   \leq C h_Y^{\dS}\Big(
   \|\dot{u}+\call_\M u\|_{ \ell^2(Y)}^2
   + h_Y^{2m-\dS-4} \|\dot{u}+\call_\M u\|_{\H^{m-2}(\M)}^2
   \Big)
   \nonumber \\
   &&\leq
   C h_Y^{\dS} \Big(
   \|\dot{u} +   \call_\M u\|_{ \ell^2(Y)}^2
   + h_Y^{2m-\dS-4}
   \big( h_Z^{-2m+4}  \|\dot{u}\|_{\H^0(\M)}^2
   + h_Z^{-2m}  \|u\|_{\H^0(\M)}^2  \big)
   \Big).
    \nonumber \\
   \label{eq:||udot-Lu||}
\end{eqnarray}
Under the denseness  constraints
\begin{equation}\label{eq:desnseness stability}
\left\{
  \begin{array}{rl}
    C h_Y^{2m- 4}  h_Z^{-2m}  < \frac14 & \mbox{if }h_Z \leq 1, \\
    C h_Y^{2m- 4}  h_Z^{-2m+4} < \frac14& \mbox{otherwise},\\
  \end{array}
\right.
\end{equation}
we integrate \eref{eq:||udot-Lu||} from $t=0$ to $T$ to finish the proof. \qed

Finally, we arrive the main result for the  convergence of the  proposed  semi-discretized solution to the  method of lines and kernel-based strong form collocation method.

\begin{theorem}\label{thm:conv}
	For some $m\geq \lfloor3+\dS/2\rfloor$, let  
$u^*\in (L^2 \cap L^\infty)(0,T;\H^m(\M))$ and $\dot{u}^*\in L^2(0,T;\H^{m-2}(\M))$
be the solution to \eref{eqPDE}.
Suppose that the surface $\M$ satisfies the assumptions in Sect.~\ref{sec:Notations and preliminaries}, kernel $\Psi_m$ satisfies \eref{kernelFour}, and two sets of quasi-uniform points $Y, \,Z\subset \M$ satisfy assumptions mentioned before and \eref{eq:desnseness stability}.
If $\alpha < \alpha^*:= h_Y^{m-\dS/2}$, further suppose that   $Y$ and  $Z$  satisfy \eref{eq:alpha*}.
Let $u_{Z,\alpha}\in \calu_{Z,\M ,\Psi_m}$ be the numerical solution defined in  \eref{eq:un} and  $\mathcal{E}_\mu {[u]} $ the residual functional defined in \eref{eq:errorfunctional}. Then, the following estimate holds
\begin{eqnarray*}
\mathcal{E}_\mu {[u_{Z,\alpha}-u^*]}
&\leq& C\Big( h_Z^{2m-4-\dS} {\int_0^T}  \|\dot{u}^*{(\cdot,\tau)}\|_{  \H^{{m-2}}(\M)} ^2 +\|u^*{(\cdot,\tau)}\|_{  \H^m(\M)}^2 {d\tau}
\\ && \qquad\qquad
+\big(  h_Z^{2m-2} + h_Y^{-2}h_Z^{2m-\dS} + h_Y^{\dS-2}\alpha^2 \big)\| u^*(\cdot,0) \|_{\H^m(\M)}^2
\Big)
\end{eqnarray*}
for $\mu=1,2$ and some constant $C$ independent of $u^*$.
\end{theorem}
\pf Consider $\mathcal{E}_\mu {[u_{Z,\alpha}-u^*]} \leq  \mathcal{E}_\mu[{u}_{Z,\alpha}-s] + \mathcal{E}_\mu[s-u^*]$;  the latter were analyzed in \eref{E(s-u*)}.
It remains to show consistency in the first term via interpolation theories.
Applying Lem.~\ref{lem:stability}, Cor.~\ref{lem:stability2}, followed by a triangle inequality  to the trial function $(u_{Z,\alpha} - s)$ yields
\begin{eqnarray*}
    \mathcal{E}_\mu {[u_{Z,\alpha}-s]}
&\leq&  C \Big(  h_Y^{\dS}  \int_{0}^{T} \|(\dot{u}_{Z,\alpha} - \dot{s}) +   \call_\M (u_{Z,\alpha} - s)\|_{ \ell^2(Y)}^2 \,{d\tau}
    + \|(u_{Z,\alpha} - s)(\cdot,0)\|_{\H^1(\M)}^2
 \Big)
 \\ &\leq& C \Big(  h_Y^{\dS}  \int_{0}^{T} \|(\dot{u}_{Z,\alpha} - \dot{s}) +   \call_\M (u_{Z,\alpha} - s)\|_{ \ell^2(Y)}^2 \,{d\tau}
 \\ &&\qquad \quad
 +  h_Y^{\dS-2} \big(
\|u_{Z,\alpha} - s\|_{ \ell^2(Y)}^2
+  \alpha^2 \|u_{Z,\alpha}\|_{\H^{m}(\M)}^2  +\| s\|_{\H^{m}(\M)}^2  \big)\Big).
\end{eqnarray*}
By Add-in and subtract-out $f=\dot{u}^*+\call_\M u^*$ and $u^*(\cdot,0)$ into the first and second norms on the right handed side, we obtain an estimate that only depends on the interpolant $s$ in \eref{eq:s} to $u^*$
\begin{eqnarray}\label{ineq:errinterp}
    \mathcal{E}_\mu {[u_{Z,\alpha}-s]}
    && \leq  2C   \Big(
     h_Y^{\dS} \int_0^T \|\dot{s} + \call_\M s - f \|_{\ell^2(Y)}^2 {d\tau}
    \nonumber \\
    && \qquad\qquad + h_Y^{\dS-2} \big(
\| s(\cdot,0) - g\|_{ \ell^2(Y)}^2
+  \alpha^2 \| s(\cdot,0) \|_{\H^{m}(\M)}^2 \big) \Big)
\label{eq:E(un-s)}
\end{eqnarray}
by the optimality \eref{eq:un} in the numerical solution.
Using the convergence estimate in \eref{eq:Izg conv} with $q=\infty$,
we have
\begin{eqnarray}
  \|\dot{s} + \call_\M s - f \|_{\ell^2(Y)}^2
  &\leq& \|\dot{s} - \dot{u}^* \|_{\ell^2(Y)}^2 + \| \call_\M s - \call_\M u^* \|_{\ell^2(Y)}^2
  \nonumber\\ &\leq&
  C h_Y^{-\dS} \Big( \|\dot{s} - \dot{u}^* \|_{L^\infty(\M)}^2 + \| \call_\M s - \call_\M u^* \|_{L^\infty(\M)}^2 \Big)
  \nonumber\\ &\leq&
  C h_Y^{-\dS}h_Z^{2m-4-\dS} \Big( \|\dot{u}^*\|_{  \H^{{m-2}}(\M)}^2
 +\|u^*\|_{  \H^m(\M)}^2 \Big).
 \label{eq:final1}
\end{eqnarray}
Next,  \eref{|IZg-g|_Y} suggests
\begin{equation}\label{eq:final2}
\| s(\cdot,0) - g\|_{ \ell^2(Y)}^2
\leq  C h_Y^{-\dS} h_Z^{2m-\dS} \| g\|_{\H^m(\M)}^2,
\end{equation}
and the orthogonality of interpolant in native space \cite[Ch.10]{Wendland-ScatDataAppr:05} suggests
\begin{equation}\label{eq:final3}
 \| s(\cdot,0) \|_{\H^{m}(\M)}^2
\leq   \| g\|_{\H^m(\M)}^2.
\end{equation}
Putting \eref{eq:final1}--\eref{eq:final3} into  \eref{eq:E(un-s)} completes the proof. \qed


\section{Fully discretized solution}
\label{sec:Fully discretized solution}

Recall the semi-discretized problem in terms of coefficients in \eref{def:uN-a} is in the from of
\begin{equation}\label{eq:model prob}
      \lambda_{Z,\alpha}(t)
    := \arginf_{\lambda,\dot\lambda\in (L^2\cap L^\infty[0,{T}])^{|{Z}|}}  \int_{0}^{T}
    \big\| \A\dot\lambda{(\tau)}+ \B\lambda{(\tau)}-\f \big \|^2_{\ell^2(\R^{|Y|})}    d{\tau}
\end{equation}
for $0\leq t \leq T$ subject to some predetermined initial condition $\lambda_{Z,\alpha}(0)$.

Suppose we numerically solve \eref{eq:model prob} at some partition $\T:=\{t_j\}$ of $[0,T]$ for $\uplambda_j$ that approximate $\lambda_{Z,\alpha}(t_j)$ by some order-$p$ scheme. Denote the corresponding fully discretized solution at any $t_j\in\T$ by
\begin{equation}\label{eq:fully discretized solution}
  U_{\T,Z,\alpha}(\cdot, t_j) :=\Psi(\cdot, Z)  \uplambda_j.
\end{equation}
Let $\zeta_j:=  \uplambda_j-\lambda_{Z,\alpha}(t_j)\in \R^{|Z|}$ and we have $ \|\zeta\|_2^2 = \calo(h_\T^{2p})$.
We can measure difference between fully and semi-discretized solutions
by
\begin{eqnarray*}
    \| U_{\T,Z,\alpha}(\cdot, t_j) - u_{Z,\alpha}(\cdot, t_j) \|_{L^2(\M)}^2
    &=&
    \zeta_j^T \Big[\int_\M \Psi(\cdot,Z)^T\Psi(\cdot,Z) \,d\sigma\Big] \zeta_j
    = \calo(h_\T^{2p}).
\end{eqnarray*}
We now derive an error estimate by a sequence of comparison
\begin{eqnarray}
    \| U_{\T,Z,\alpha}(\cdot, t_j) - u^*(\cdot, t_j)
    \|_{\ell^\infty(\T;L^2(\M))}^2
     &\leq&
    \| U_{\T,Z,\alpha}(\cdot, t_j) - u_{Z,\alpha}(\cdot, t_j) \|_{\ell^\infty(\T;L^2(\M))}^2
    \nonumber\\&&\qquad\qquad
    +  \| u_{Z,\alpha}(\cdot, t_j) - u^*(\cdot, t_j) \|_{\ell^\infty(\T;L^2(\M))}^2
    \nonumber\\
    &\leq &
    \| U_{\T,Z,\alpha}(\cdot, t_j) - u_{Z,\alpha}(\cdot, t_j) \|_{\ell^\infty(\T;L^2(\M))}^2
    \nonumber\\&&\qquad\qquad
    + \esssup_{0\leq {\tau}\leq T}\| u_{Z,\alpha}{(\cdot,\tau)} - u^*(\cdot, {\tau}) \|_{ H^1(\M)}^2
    \nonumber\\
    &\leq &
    \calo(h_\T^{2p}) + \mathcal{E}_\mu {[u_{Z,\alpha}-u^*]},
    \label{eq:full error}
\end{eqnarray}
allowing Theorem~\ref{thm:conv} to be applied to the fully discretized solution.

\subsection{Temporal discretization}

Despite being linear, the second order Euler-Lagrange ODE for \eref{eq:model prob} involves multiple products of kernel matrices, which is ill-conditioned. We propose two  variants of algorithms for solving \eref{eq:model prob}: (i) difference equation and (ii) ODE approach.
We focus on uniform time discretization at
\begin{equation}\label{eq:T}
\T = \{t_j \}_{j=0}^{|\T|}:= \{ j h_\T\}_{j=0}^{|\T|}.
\end{equation}
Further suppose that the integrand of \eref{eq:model prob}, i.e.,  squares of norms of PDE residuals, is of class $C^{p}[0,T]$. This assumption imposes temporal smoothness requirements on $f$ in \eref{eqPDE} and hence to $u^*$.
Then, there exists a  order-$p$ backward finite difference scheme $\Dh$, see \cite[Tab. 3]{Fornberg-Genefinidiffform:88} for coefficients, such that
\[
    \dot\uplambda(t_j)= \Dh \uplambda_j + \calo( h_\T^{p} )
    =\frac{1}{h_\T}\Big( \gamma_0\uplambda_j +\sum_{k=-p}^{-1}\gamma_k \uplambda_{k+j} \Big)+\calo( h_\T^{p} ),
\]
and
an order-$p$ numerical quadrature rule $\{w_j\}$ on $\T$, say composite Newton-Cotes rules, so that, given $\lambda_{Z,\alpha}(0)$, we can approximate the set of solutions \eref{eq:model prob} at $\T$ by

\begin{eqnarray}
\lambda_{Z,\alpha}(\T)
    &:=& \left\{ \arginf_{\lambda,\dot\lambda\in (L^2\cap L^\infty[0,{T}])^{|{Z}|}}  \int_{0}^{T}
    \big\| \A\dot\lambda{(\tau)}+ \B\lambda{(\tau)}-\f \big \|^2_{\ell^2(\R^{|Y|})}    {d\tau}  \right\}_{\big|\T}
    \nonumber\\
    &=& \left\{ \arginf_{\substack{ \uplambda_j \in \R^{|Z|} \\ 1\leq j\leq |\T|}}
    \sum_{t_j\in\T} w_j
    \big\| \A\Dh \uplambda_j  + \B\uplambda_j-\f_j  \big \|^2_{\ell^2(\R^{|Y|})}      \right\}
    +\calo( h_\T^p ).  \label{def:discretize in time 2}
\end{eqnarray}
Finally, we arrive at a sum-of-squares problem in \eref{def:discretize in time 2}, in which we found standard least-squares problems with respect to $\uplambda_j$. Explicitly, we write the solution to \eref{def:discretize in time 2}  by the following recursive formula:
\begin{subequations}
	\begin{eqnarray}
		&&\qquad\qquad\lambda_{Z,\alpha}(\T)
		=
		\left\{ \arginf_{  \uplambda_j \in \R^{|Z|}  }
		\big\| \A\Dh \uplambda_j  + \B\uplambda_j-\f_j  \big \|^2_{\ell^2(\R^{|Y|})}    \right\}_{ j=1}^{|\T|}
		\!\!+\calo( h_\T^p )
		\label{eq:LS lambda0}
		\\&&\qquad\qquad=
		\left\{ \arginf_{  \uplambda_j \in \R^{|Z|}  }
		\left\|\Big(\frac{ {\gamma_0}}{h_\T}\A + \B\Big) \uplambda_j + \frac{1}{h_\T}
		\A { \Big( \sum_{k=-p}^{-1}\gamma_j \uplambda_{k+j}\Big) } -\f_j\right\|^2_{\ell^2(\R^{|Y|})} \right\}_{ j=1}^{|\T|}
		+\calo( h_\T^p )
		\label{eq:LS lambda1}
		\\&&\qquad\qquad=
		\left\{
		\Big(\frac{ {\gamma_0}}{h_\T}\A + \B\Big)^\dagger
		\left(  \f_j - \frac{1}{h_\T}\A { \Big( \sum_{k=-p}^{-1}\gamma_j \uplambda_{k+j}\Big) }  \right)
		\right\}_{ j=1}^{|\T|}
		+\calo( h_\T^p ),
		\label{eq:LS lambda}
	\end{eqnarray}
\end{subequations}
subject to initial condition $\uplambda_0 = \lambda_{Z,\alpha}(0) \in \R^{|Z|}$.
We summarize the above findings.

\begin{theorem}\label{thm:conv 2}
  Suppose that the assumptions in Thm.~\ref{thm:conv} hold. For some integer $p\geq 1$,
  we further suppose that
 $u^*\in C^{p+1 }(0,T;\H^m(\M))$.
For some $p$-th order backward finite difference scheme on uniform time grid $\T$ in \eref{eq:T} for the first derivative, let its coefficients be denoted as $\gamma_{-p},\ldots,\gamma_0\in\R$.
   Let $\uplambda_0 = \lambda_{Z,\alpha}(0) \in \R^{|Z|}$ as in \eref{def:uN-b} and
  \begin{equation}\label{eq:lam_j}
      \uplambda_j=
  \Big(\frac{ {\gamma_0}}{h_\T}\Psi_m(Y,Z) + [\call_\M\Psi_m](Y,Z)\Big)^\dagger
    \left( f(Y,t_j) - \frac{1}{h_\T}  \Psi_m(Y,Z)   { \Big( \sum_{k=-p}^{-1}\gamma_j \uplambda_{k+j}\Big) } \right),
  \end{equation}
  for $1\leq j \leq|\T|$,
  be the coefficients of the fully discretized solution $U_{\T,Z,\alpha}(\cdot,t_j)$ in \eref{eq:fully discretized solution}. Then, the error estimate in Thm.~\ref{thm:conv} for semi-discretized solution extends to the fully  discretized solution with an additional $p$-th order temporal error term, i.e.,
  \[
      \| U_{\T,Z,\alpha}(\cdot, t_j) - u^*(\cdot, t_j)
    \|_{\ell^\infty(\T;L^2(\M))}^2
     \leq
    \calo(h_\T^{2p}) + \mathcal{E}_\mu {[u_{Z,\alpha}-u^*]}
    \]
     holds.
\end{theorem}

Note that the difference equations in \eref{eq:LS lambda} and  \eref{eq:lam_j} were derived without any ODE in the process. To benefit from the vast library of ODE solvers with adaptive time stepping, we conclude our theoretical work with an equivalent ODE.

\begin{corollary}\label{cor ode}
  Thm.~\ref{thm:conv 2} remains to hold if the difference solution \eref{eq:lam_j} were replaced by some $p$-th order approximation to the solution of the following $|Z|\times |Z|$ systems of ODE
  \begin{equation}\label{eq: final ode}
    \dot\uplambda(t) = \Psi_m(Y,Z)^\dagger \Big( f(Y,t) -  [\call_\M\Psi_m](Y,Z) \uplambda(t) \Big)
  \end{equation}
  for $\uplambda\in \big(C^{p+1}[0,T]\big)^{|Z|}$.
\end{corollary}
\pf Using notations in \eref{eq:LS lambda}, the least-squares solution \eref{eq:lam_j} solves the normal equation
\[
   \Big(\frac{ {\gamma_0}}{h_\T}\A + \B\Big)^T
   \left( \Big(\frac{ {\gamma_0}}{h_\T}\A + \B\Big)\uplambda_j
    + \frac{1}{h_\T}\A { \Big( \sum_{k=-p}^{-1}\gamma_j \uplambda_{k+j}\Big) }
    -\f_j    \right) = 0
\]
for any give partition $\T$. By reverse-manipulation from \eref{eq:LS lambda1} back to \eref{eq:LS lambda0}, we can rewrite the normal equation as below and take limit of $h_\T = T/|\T| \to 0$:
\[
   \lim_{h_\T\to0} \Big({ {\gamma_0}}\A + {h_\T}\B\Big)^T
   \lim_{h_\T\to0}  \Big( \A\Dh \uplambda_j + \B \uplambda_j   -\f_j   \Big) = 0,
   \quad 1\leq j \leq |\T|,
\]
to obtain
\[
   \A ^T
   \Big( \A\dot\uplambda(t) + \B \uplambda(t)  -\f(t)   \Big) = 0   ,
   \quad 0<t\leq T,
\]
which is the normal equation of \eref{eq: final ode}. \qed


\section{Numerical examples}\label{sec:Numerical examples} This section contains examples with self-explanatory title designed for various numerical verifications and aims. We use $y=(x_1,x_2,x_3)\in\M$ to define functions throughout the section.

\subsection*{Example 1: Comparing with a meshless Galerkin method}

\begin{table}
		\centering
	\caption{{Exmp.~1: Relative $L^2(\M)$-errors at $T=1$,  estimated order of convergence (eoc) of solution \eref{eq:lam_j} with $p=2$ and $m=4$, and a graphical comparison between  error of a meshless Galerkin method\cite[Tab.~2]{Kunemund+NarcowichETAL-highmeshGalemeth:19} using  $3721$ thin-plate splines  with various numbers of quadrature points (in {\color[RGB]{51,51,255} blue}, {\color[RGB]{255,0,0} red}, and {\color[RGB]{0,0,0} black} as in \cite[Fig.~1 (Right)]{Kunemund+NarcowichETAL-highmeshGalemeth:19}) and the reported RBF-MOL error (in {\color[RGB]{102,153,51} green}), which are not distinguishable to the eyes.}}\label{TableRelErrSphere}
\begin{minipage}{0.7\linewidth}
	\begin{tabular}{lllcllcll}
			\hline
		\multirow{2}{*}{$|Y|$}   & \multirow{2}{*}{$h_\tau$} & \multicolumn{2}{l}{$|Z|=961$} & &\multicolumn{2}{l}{$|Z|=3721$} \\\cline{3-4}
	\cline{6-7}	&               &         Rel. error           & eoc    &   & Rel. error             & eoc       \\
			\hline\hline
		\multirow{4}{*}{\rotatebox{90}{1,153}}
        & 0.06                             & 1.198\blue{428}E-4      & -    &     &   &       \\
		& 0.04                             & 5.602\blue{314}E-5      & 1.88    &  &  &      \\
		& 0.02                             & 1.2504\blue{87}E-5      & 2.16  &    &  &      \\
		& 0.01                             & 2.9277\blue{79}E-6      & 2.09   &   &  &      \\	\hline
		\multirow{4}{*}{\rotatebox{90}{4,465}}
        & 0.06                             & 1.198337E-4             & -     &             & 1.198\blue{429}E-4  & -         \\
		& 0.04                             & 5.60197\blue{6}E-5      & 1.88      &         & 5.602\blue{329}E-5  & 1.88      \\
		& 0.02                             & 1.250414E-5             & 2.16     &          & 1.250\blue{503}E-5  & 2.16      \\
		& 0.01                             & 2.92770\blue{9}E-6      & 2.09      &         & 2.927\blue{935}E-6  & 2.09      \\	\hline  \hline
		\multirow{4}{*}{\rotatebox{90}{23,042}}
        & 0.06                             & 1.198337E-4             & -     &    & 1.198337E-4  & -         \\
		& 0.04                             & 5.60197\blue{6}E-5      & 1.88    &  & 5.601977E-5  & 1.88      \\
		& 0.02                             & 1.250414E-5             & 2.16  &    & 1.250414E-5  & 2.16      \\
		& 0.01                             & 2.92770\blue{9}E-6      & 2.09   &   & 2.927710E-6  & 2.09      \\	\hline
		\multirow{4}{*}{\rotatebox{90}{40,962}}
        & 0.06                             & 1.198337E-4             & -        & & 1.198337E-4  & -         \\
		& 0.04                             & 5.60197\blue{6}E-5      & 1.88     & & 5.601977E-5  & 1.88      \\
		& 0.02                             & 1.25041\blue{3}E-5      & 2.16     & & 1.250414E-5  & 2.16      \\
		& 0.01                             & 2.9277\blue{00}E-6      & 2.09     & & 2.927710E-6  & 2.09      \\	\hline
		\multirow{4}{*}{\rotatebox{90}{92,162}}
        & 0.06                             & 1.198337E-4             & -         && 1.198337E-4  & -         \\	
		& 0.04                             & 5.60197\blue{6}E-5      & 1.88      && 5.601977E-5  & 1.88      \\
		& 0.02                             & 1.25041\blue{3}E-5      & 2.16      && 1.250414E-5  & 2.16      \\
		& 0.01                             & 2.92770\blue{5}E-6      & 2.09      && 2.927710E-6  & 2.09      \\	\hline
		\multirow{4}{*}{\rotatebox{90}{256,002}}
        & 0.06                             & 1.198337E-4             & -         &&  1.198337E-4  & -         \\
		& 0.04                             & 5.60197\blue{6}E-5      & 1.88      &&  5.601977E-5  & 1.88      \\
		& 0.02                             & 1.25041\blue{3}E-5      & 2.16      &&  1.250414E-5  & 2.16      \\
		& 0.01                             & 2.92770\blue{5}E-6      & 2.09      &&  2.927710E-6  & 2.09     \\    \hline
	\end{tabular}
\end{minipage}
\begin{minipage}{0.28\linewidth}
	\begin{overpic}[width=\textwidth,trim=0 0 0 0, clip=true,tics=10]{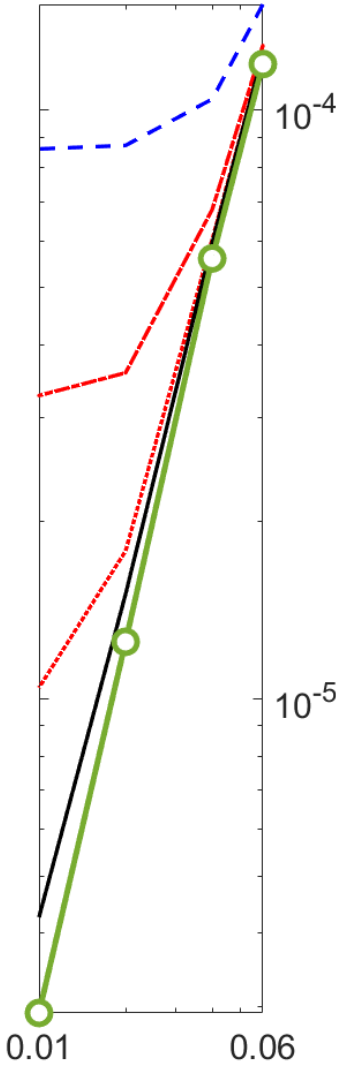}
\put(4.2,87){\scriptsize \rotatebox{90}{\color[RGB]{51,51,255} 23,042}}
\put(4.2,64){\scriptsize \rotatebox{90}{\color[RGB]{255,0,0} 40,962}}
\put(4.2,40){\scriptsize \rotatebox{90}{\color[RGB]{255,0,0} 92,162}}
\put(4.2,23){\scriptsize \rotatebox{90}{256,002}}
\put(7,10){{\color[RGB]{102,153,51} RBF-MOL}}
\end{overpic}
\end{minipage}
\end{table}

Our first aim is to compare the proposed method with a weak formulation in \cite{Kunemund+NarcowichETAL-highmeshGalemeth:19}. We consider an example there: $\dot u - 0.1\Delta_\M {u} + 3 u = f$
on the unit sphere with  $u^*(y,t) =\exp( x_1+1/(1+t) )$ for $t\in[0,1]$.
We solve the problem by $m=4\geq \lfloor3+\dS/2\rfloor$ Sobolev kernel and the difference equations in \eref{eq:lam_j} 
with $p=2$.
{
We obtain numerical approximations by the difference equations in \eref{eq:lam_j}.
We use the same time stepping $0.01\leq h_\T \leq 0.06$ and  number of trial centers $|Z|\in\{961,\,3721\}$ as in \cite{Kunemund+NarcowichETAL-highmeshGalemeth:19}. We use collocation/quadrature point sets of size $|Y|=\{1153,\,4465,\,23042,\,\ldots,256002\}$, in which the first two tested $|Y|$ where selected by $120\%|X|$ and the rest were used in \cite{Kunemund+NarcowichETAL-highmeshGalemeth:19}.

In Tab.~\ref{TableRelErrSphere}, we report the relative discrete $L^2(\M)$-error  of the proposed method. The errors of our methods are obviously limited by the time discretization. That is,  increasing $|Z|$ does not improve accuracy for all tested $|Y|\geq 120\%|Z|$.
For $m=4$, one can indeed use the Lagrange  interpolation setup \cite{Wendland2010A} with $X=Y$ in method of lines (MOL) and obtain errors similar in magnitudes as in $|Y| = 120\%|Z|$. Oversampling become a must when $m$ and/or $|Z|$ become large. Due to page limitation, we refer readers to \cite{Chen+Ling-OverneceRBF-meth:22} for detailed numerical studies.
Errors of all tested cases agree with that of our finest resolution up to  4 or more significant figures.
The observed estimated orders of convergence (eoc) support the theoretical temporal convergence rate in Thm.~\ref{thm:conv 2}.

Since the meshless Galerkin method \cite{Kunemund+NarcowichETAL-highmeshGalemeth:19},  which uses  Crank-Nicolson (CN) and   thin-plate spline kernels, did not contain errors for larger  $m$,  we cannot fairly compare the two methods, i.e., their numerical denseness requirement of quadrature points.
Yet, it is safe to say that these methods are comparable in terms of accuracy provided that $|Y|$ is sufficiently large with respect to $|Z|$.
}

\subsection*{Example 2: Nonconstant diffusive tensor and kernels' smoothness}

\begin{figure}
\begin{minipage}[t]{0.45\textwidth}
	\centering
	\begin{overpic}[width=\textwidth,trim=0 0 0 0, clip=true,tics=10]{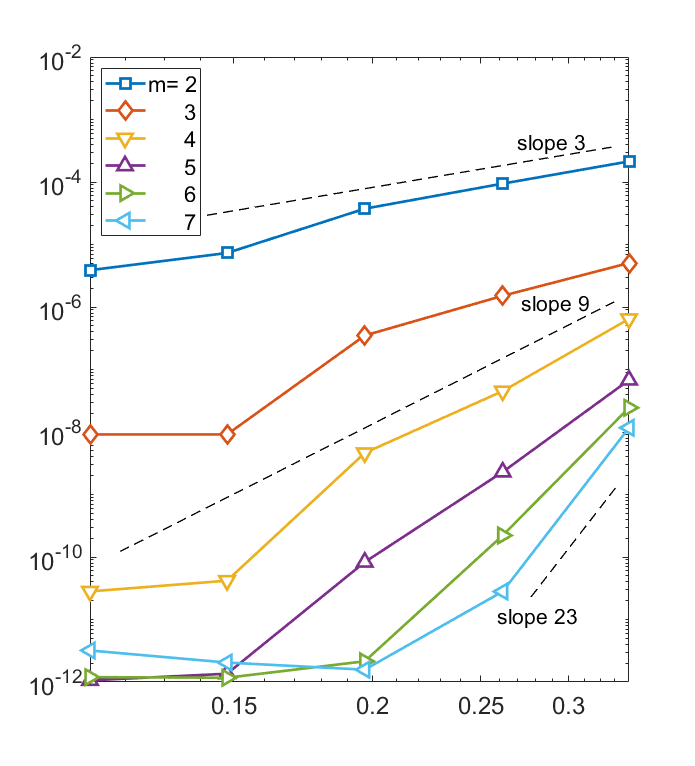}
\put(45,2){\scriptsize${h_Z}$}
\put(0,40){\scriptsize \rotatebox{90}{Rel. $L^2(\M)$-error}}
\end{overpic}
	\caption{Exmp.~2: Spatial convergence profiles with respect to kernel smoothness.}\label{Fig:NonIdentityTensor}
\end{minipage}
\begin{minipage}[t]{0.45\textwidth}
	\centering
	\begin{overpic}[width=\textwidth,trim=0 0 0 0, clip=true,tics=10]{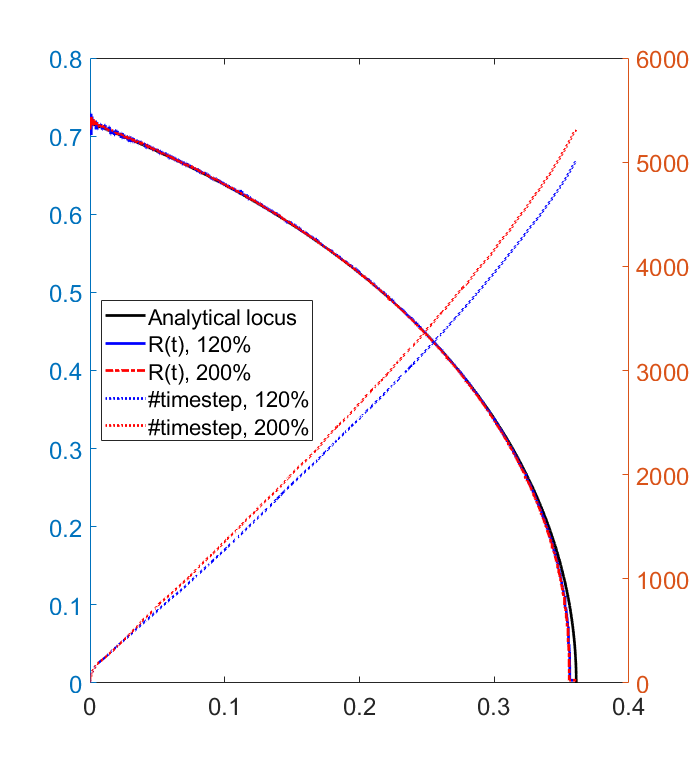}
\put(45,2){\scriptsize$t$}
\put(0,40){\scriptsize \rotatebox{90}{Radius, R(t)}}
\put(93,58){\scriptsize \rotatebox{270}{\#timestep}}
\put(40,16){\color{black}\frame{\includegraphics[scale=.35]{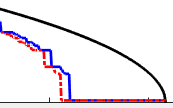}}}
\put(65,16){\color{black}\vector(2,-1){8}}
\end{overpic}
	\caption{Exmp.~3: Radius the spherical cap and number of time stepping used by ODE45.}\label{Fig:AllenCahn_cap}
\end{minipage}
\end{figure}

Using the projection matrix $P$ in \eref{def:P(y)}, let  $A(y)=P(y)\,\text{Diag}[x_1^2+1,1,1]\in\R^{3\times3}$ be the diffusion tensor.
We solve the diffusion equation \eref{eqPDE}
on the unit sphere {to $T=1$}. We use $u^*$ in Exmp.~1 as exact solution. This time around, we solve the PDE via ODE \eref{eq: final ode} in Cor.~\ref{cor ode}.
We perform a short-term temporal integration by the implicit Crank-Nicolson method presented in \cite{Kunemund+NarcowichETAL-highmeshGalemeth:19}
with a sufficient small time-step size $h_\T$, {which is $1\rm{E}-6$ in this example,}
in order to study the spatial convergence behaviour of the proposed method with respect to various kernel's smoothness $2\leq m \leq 7$.
{Once again, we see that it is numerically sufficient to use linear ratio of oversampling to observe spatial convergence; using more collocation points in $Y$ to satisfy the theoretical denesness/regularization requirment in \eref{eq:desnseness stability} will (very) slightly improve accuray (as in Exmp.~1).}
Relative $L^2(\M)$-error by using $|Z|=100$ to 1000 trial centers and 120\% oversampling were shown in Fig.~\ref{Fig:NonIdentityTensor}.
Generally speaking, larger   $m$  yields   faster  (inital) eoc before error stagnation.
Although $m=2$ and 3 are not covered by Thm.~\ref{thm:conv}, convergence is still observed.
Also, it is not uncommon to observe faster than theoretically predicted eoc when $u^*$ and $\M$ are both of high order of smoothness, see numerical experiment in \cite[Exmp. 2]{Cheung+Ling-Kernembemethconv:18}.

\subsection*{Example 3: Simulating Allen-Cahn equations}
In this numerical experiment, we consider the Allen-Cahn equation
$
\dot{u}= \Delta_\M u + \frac{1}{\varepsilon^2}u(1-u^2),
$
which is a reaction diffusion equation that models phase separation of two fluids\cite{ALLEN19791085} with $\varepsilon>0$ being the width of the diffusion interface between two fluids.
We set up the test problem on the unit sphere as in \cite{Kunemund+NarcowichETAL-highmeshGalemeth:19,Choi-AC_sphere:15} {with $T:=-\frac{1}{2}\log(1-R_0^2)$ for some $R_0>0$}.
We define an initial condition $u(0)$ to be  $+1$ within radius $R_0$ to the northpole and $-1$ otherwise.
Then,
$R(t)=\sqrt{ 1 - ( 1- R_0^2)e^{2t}  }$
is the radius  {of}  the shrinking spherical cap of the solution for $0\leq t\leq T$.

\begin{figure}
	\centering
	\begin{overpic}[width=0.19\textwidth,trim=50 50 80 45,
		clip=true,tics=10]{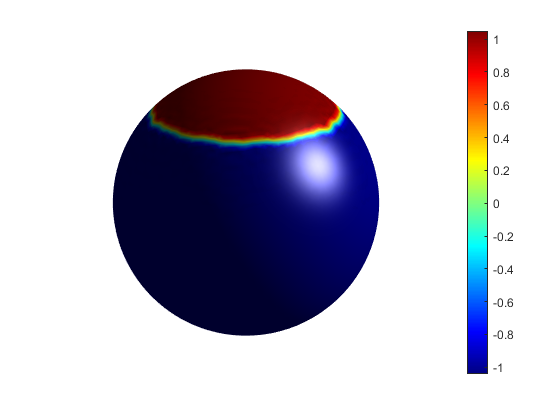}
		\put(28,80){\scriptsize$t=5E-4$}\end{overpic}
	\begin{overpic}[width=0.19\textwidth,trim=50 50 80 45, clip=true,tics=10]{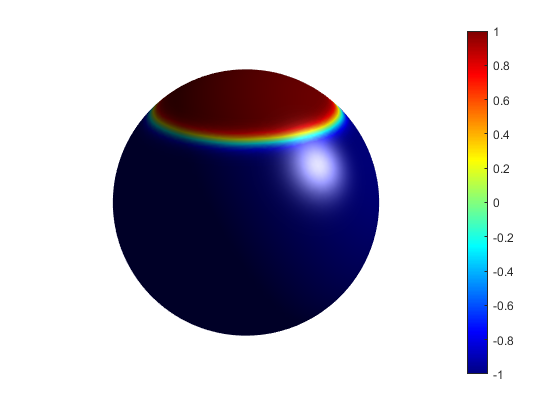}
		\put(35,80){\scriptsize$0.01$}\end{overpic}
	\begin{overpic}[width=0.19\textwidth,trim=50 50 80 45, clip=true,tics=10]{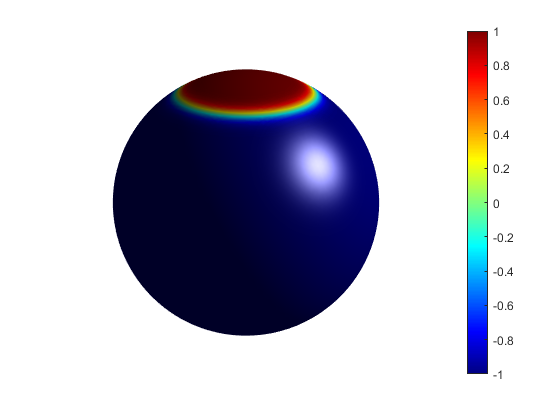}
		\put(35,80){\scriptsize$0.20$}\end{overpic}
	\begin{overpic}[width=0.19\textwidth,trim=50 50 80 45, clip=true,tics=10]{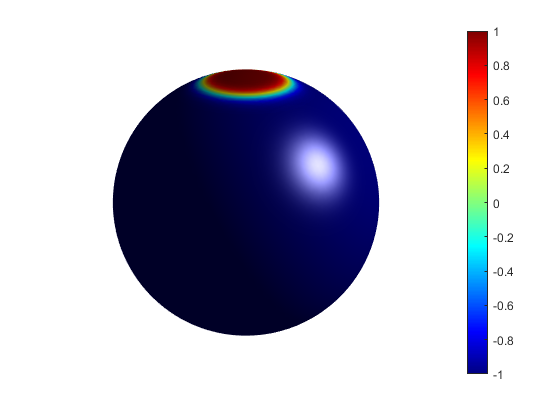}
		\put(35,80){\scriptsize$0.30$}\end{overpic}
	\begin{overpic}[width=0.19\textwidth,trim=50 50 80 45, clip=true,tics=10]{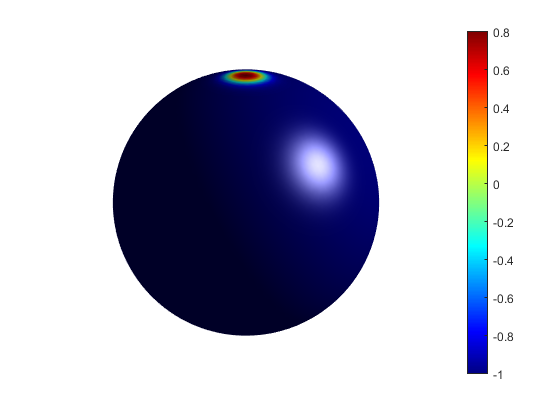}
		\put(35,80){\scriptsize$0.35$}\end{overpic}	
	\caption{Exmp.~3: Snapshots of the shrinking cap towards the north pole.}\label{Fig:AllenCahn}
\end{figure}

Using $R_0=0.717$, $\varepsilon=0.05$, and $|Z|=3721$ as in \cite{Kunemund+NarcowichETAL-highmeshGalemeth:19} for another round of comparison with the Galerkin method, we solve the semi-discrete equation (\ref{eq: final ode}) by the MATLAB built-in ODE45 solver. Firstly, we compute the econ-QR factorization of $\Psi_m(Y,Z)=QR$ where $Q$ and $R$ are of size $|Y|\times|Z|$ and $|Z|\times|Z|$. Replacing  $\Psi_m(Y,Z)^\dagger$ in \eref{eq: final ode} by $Q^T$, we call ODE45 with mass matrix $R$.
We test oversampling ratios $120\%$ and $200\%$ with $|Y|=4465$ and $7442$;
Fig.~\ref{Fig:AllenCahn_cap}~(Left-$y$) shows the analytical and numerical radius of the spherical cap, which should be compared with \cite[Fig. 3]{Kunemund+NarcowichETAL-highmeshGalemeth:19}. First of all, both kernel-based methods of lines show good accuracy for a long time. If we want to be really picky,
\begin{itemize}
  \item the proposed method shows small oscillations in radius at small $t$ due to fewer numbers of collocation points, but
  \item the Galerkin solutions with $|Y|=40962$ and $92162$ have visually different (without zoom-in) final cap-vanishing times.
\end{itemize}
Fig.~\ref{Fig:AllenCahn_cap}~(Right-$y$) shows the numbers of time steps required and
Fig.~\ref{Fig:AllenCahn} shows a few snapshots of the shrinking spherical cap.

%
%
%

\begin{figure}
	\centering
\begin{overpic}[width=0.19\textwidth,trim=50 50 40 40, clip=true,tics=10]{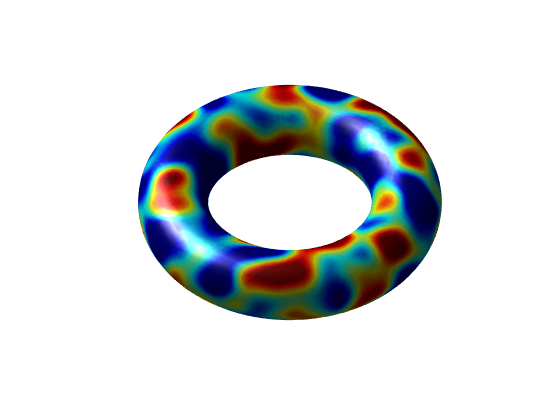}
\put(-10,10){\scriptsize \rotatebox{90}{$|Z_1|=864$}}
    \put(32,75){\scriptsize$t=0.006$}\end{overpic}
\begin{overpic}[width=0.19\textwidth,trim=50 50 40 40, clip=true,tics=10]{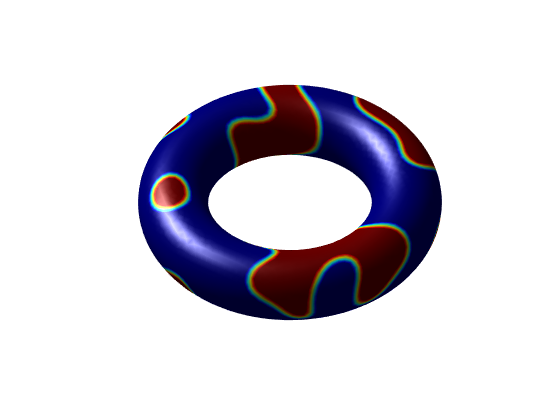}
    \put(40,75){\scriptsize$0.02$}\end{overpic}
\begin{overpic}[width=0.19\textwidth,trim=50 50 40 40, clip=true,tics=10]{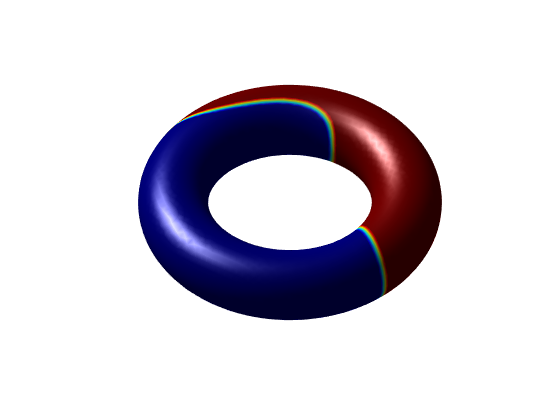}
    \put(42,75){\scriptsize$0.5$}\end{overpic}	
\begin{overpic}[width=0.19\textwidth,trim=50 50 40 40, clip=true,tics=10]{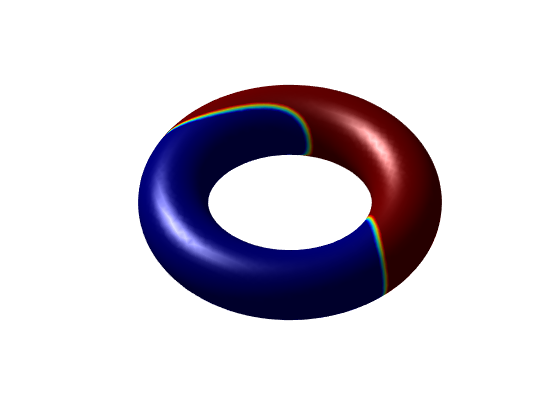}
    \put(42,75){\scriptsize$1.0$}\end{overpic}
\begin{overpic}[width=0.19\textwidth,trim=50 50 40 40, clip=true,tics=10]{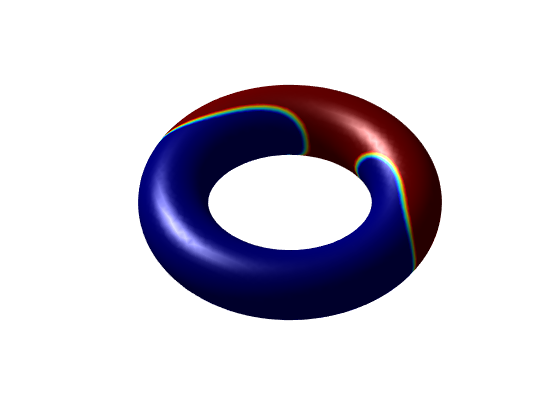}
    \put(42,75){\scriptsize$4.0$}\end{overpic}
\\
\begin{overpic}[width=0.19\textwidth,trim=50 50 40 40, clip=true,tics=10]{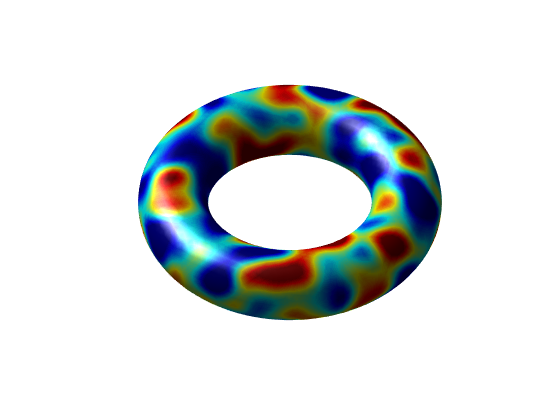}
\put(-10,10){\scriptsize \rotatebox{90}{$|Z_2|=1696$}}
     \end{overpic}
\begin{overpic}[width=0.19\textwidth,trim=50 50 40 40, clip=true,tics=10]{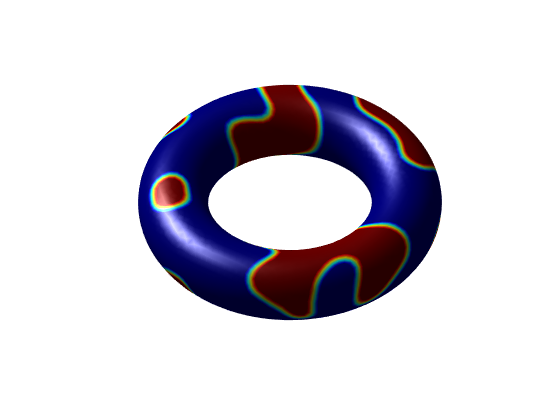}
     \end{overpic}
\begin{overpic}[width=0.19\textwidth,trim=50 50 40 40, clip=true,tics=10]{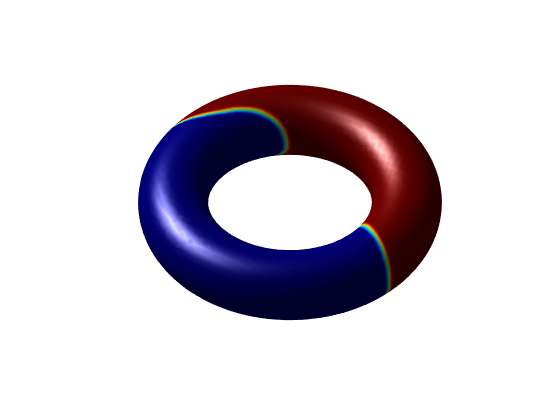}
     \end{overpic}	
\begin{overpic}[width=0.19\textwidth,trim=50 50 40 40, clip=true,tics=10]{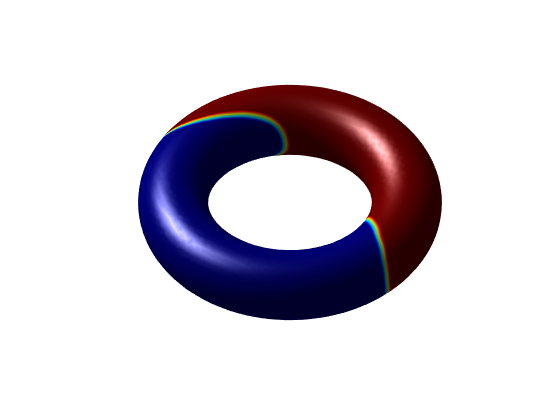}
     \end{overpic}
\begin{overpic}[width=0.19\textwidth,trim=50 50 40 40, clip=true,tics=10]{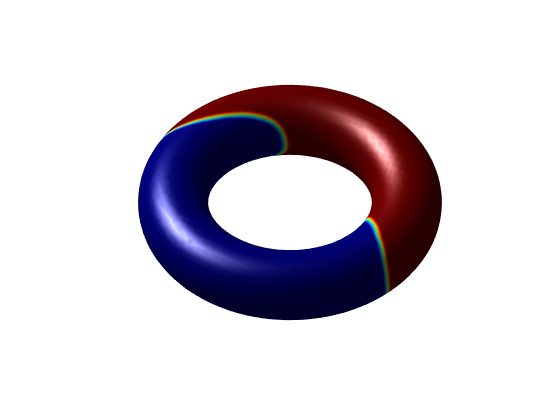}
     \end{overpic}
	\caption{Exmp.~3:  Snapshots of Allen-Cahn solutions based on the same random initial condition and collocation points with size $|Y|=3856$, but different sets of trial centers.}\label{Fig:AllenCahn:separtion}
\end{figure}

We end this example with an Allen-Cahn phase separation simulation on a torus\footnote{Torus: $x^2 + y^2 + z^2 + 1^2-(1/3)^2)^2-4(x^2 + y^2) = 0$}.
Using the generator in \cite{Persson+Strang-SimpMeshGeneMATL:04} with $h_{Z_1}=0.1333$ and $h_{Z_2}=0.1$, and  $h_X=0.066$, we obtain point sets of sizes $|Z_1|= 864$, $|Z_2|= 1696$, and $|X|=3856$ to repeat the same calculations.
Random initial values between $[-0.5, 0.5]$ were assigned to $X$ and we present in Fig.~\ref{Fig:AllenCahn:separtion} a selected initial condition that leads to steady state  solutions with non-vanishing diffuse interface. Comparing results in the eyeball norm, two trial solutions look obviously different at $t=0.5$. We also note that the solution of $Z_2$ with more trial centers arrives steady state earlier (c.f. $t=1$ and $4$).

\section{Conclusions}
We theoretically deduce a method of lines based on a discrete least-squares approximation to the initial condition and strong-form collocation with sampling for diffusion equations on smooth surfaces.
Convergence estimate of the (spatially) semi-discretized solution is given in Thm.~\ref{thm:conv}.
The fully discretized problem is a difference equation (see Thm.~\ref{thm:conv 2}) that were connected to a system of ODE (see Cor.~\ref{cor ode}).
{
Numerical  examples were provided to  compare the proposed method with a meshless Galerkin method.}
We remind readers that our strong-form theories require  solutions with higher smoothness, i.e., $m\geq \lfloor3+\dS/2\rfloor$. For PDEs with lower regularity, meshless Galerkin methods should be the theoretically sound method of choice.

\section*{Acknowledgements}
This work was supported by the General Research Fund (GRF No. 12301917, 12303818, 12301419) of Hong Kong Research Grant Council, National Natural Science Foundation (Grant No. 12001261)  {and Jiangxi Provincial Natural Science Foundation (Grant No. 20212BAB211020)}.

\bibliographystyle{siamplain} 
\bibliography{SINUM144436R2G}

\end{document}